\documentclass[preprint]{imsart}
\RequirePackage[OT1]{fontenc}
\RequirePackage{amsthm,amsmath}
\usepackage{amsfonts}
\usepackage{amssymb}
\usepackage{lineno}
\RequirePackage[numbers]{natbib}
\RequirePackage[colorlinks,citecolor=blue,urlcolor=blue]{hyperref}

\startlocaldefs

\numberwithin{equation}{section}
\theoremstyle{plain}

\newtheorem{theorem}{Theorem}[section]
\newtheorem{lemma}[theorem]{Lemma}

\newtheorem{corollary}[theorem]{Corollary}
\newtheorem*{condition}{Condition}
\newtheorem*{conditionA}{Condition A}
\newtheorem*{conditionB}{Condition B}

\newtheorem{remark}{Remark}
\newtheorem{assumption}{A.}[section]

\usepackage{tikz}
\newcommand*\circled[1]{\tikz[baseline=(char.base)]{
		\node[shape=circle,draw,inner sep=1pt] (char) {#1};}}
\usepackage{enumitem}

\endlocaldefs
\begin{document}

\begin{frontmatter}
	%\title{A Z-estimation system for model development with two-phase sampling}
	\title {Z-estimation system: a modular approach to asymptotic analysis}
	
	\begin{aug}
		\author{Jie Kate Hu}
		
		\affiliation{Department of Biostatistics, Harvard T.H. Chan School of Public Health }

 \today
% \thankstext{t1}{The author thanks Norm Breslow, Jon Wellner, and Gary Chan for valuable discussions that led to this paper.} 

	%%	\ead[label=e1]{}
		%\ead[label=u1,url]{}}
		\end{aug}

	%\runauthor{}
	
	%\affiliation{ }

		%\address{Address of the Third author\\
		%	\printead{e3}\\
		%	\printead{u1}}

\begin{abstract}
Asymptotic analysis for related inference problems often involves similar steps and proofs. These intermediate results could be shared across problems if each of them is made self-contained and easily identified. However, asymptotic analysis using Taylor expansions is limited for result borrowing because it is a step-to-step procedural approach. This article introduces EEsy, a modular system for estimating finite and infinitely dimensional parameters in related inference problems. It is based on the infinite-dimensional Z-estimation theorem, Donsker and Glivenko-Cantelli preservation theorems, and weight calibration techniques. This article identifies the systematic nature of these tools and consolidates them into one system containing several modules, which can be built, shared, and extended in a modular manner. This change to the structure of method development allows related methods to be developed in parallel and complex problems to be solved collaboratively, expediting the development of new analytical methods. This article considers four related inference problems --- estimating parameters with random sampling, two-phase sampling, auxiliary information incorporation, and model misspecification. We illustrate this modular approach by systematically developing 9 parameter estimators and 18 variance estimators for the four related inference problems regarding semi-parametric additive hazards models. Simulation studies show the obtained asymptotic results for these 27 estimators are valid. In the end, I describe how this system can simplify the use of empirical process theory, a powerful but challenging tool to be adopted by the broad community of methods developers. 
  %I compare this system to  Lego system, suggesting potential new benefits of this system for future development.
 I discuss challenges and the extension of this system to other inference problems. 

\end{abstract}

\end{frontmatter}
\section{Introduction}

%Modern sensors, mobile,  and web technology allow us to monitor conditions, actions, and outcomes of almost every activity, generating a vast amount of information. This information has motivated numerous new questions in science and industry. Consequently, the demand for statistical methods has exploded. However, the development process of statistical methods has not changed much over the past decades--- laboriously developing one method at a time for a specific sampling frame. There is a shortage of data scientists and statisticians in almost every industry and field of science. Can our current art of method development catch up to the rapidly emerging questions in science and industry? This article addresses these challenge from new perspectives. Instead of one method at a time, we develop an estimation system in which related problems can be  solved systemically and simultaneously. This system is also  a development platform,  on which intermediate results capsulized in modules could be easily  transported from one method to another, expediting new  method development.

Statistical methods have grown increasingly complex as the field evolves, owing to the interplay of escalation in  complexity of models, data structures, and techniques to handle  unideal data. Consequently, the  asymptotic studies   of  these methods  have  also  become extensive and intricate, and sometimes difficult to follow and verify, potentially impeding their transparency and reliability. One reason behind these problems could be the prevalence of the Taylor expansion approach for asymptotic studies, which follows a linear developmental process. It is suitable for simple inference problems, but when problems become complex, it is difficult to isolate the intermediate results in this approach for better understanding and validation. Interestingly, upon closer examination, many related inference problems studied by Taylor expansions repeat similar steps and proofs. Without altering the asymptotic analysis paradigm, the perpetual growth in method complexity renders the corresponding theoretical investigations and peer-review process unsustainable. This article addresses this issue by  switching to a modular  approach, built upon Huber's alternative asymptotic analysis method in his 1967 landmark paper\citep{huber1967}.The whole  estimation process is decomposed into  small modules, conquered individually and assembled in the end for desirable results. Small tasks in small components are easy to follow and check. This transparency, as a result, can improve the quality of methodology research. Because intermediate results are capsulized in modules,  they are easy to be identified, transported from one method to another, and  adapted for new problems , expediting new  method development.

%Two-phase sampling design was commonly used in surveys and had been gradually adopted in epidemiological studies due to its cost-effectiveness. It uses outcomes and/or inexpensive covariates obtained for each subject at Phase I to determine their following sub-sampling probabilities, so that certain expensive and difficult to be obtained covariates are only measured for the most informative subjects in the Phase II subsamples. 
Two-phase sampling design collects data in a sequential way--- a large phase I random sample  of inexpensive variables is first drawn  and stratified. Different amounts of  subsamples from different strata are then selected for expensive variable measurements.   By this means,  the cost of a study is allocated  more efficiently to the more informative subjects.  For its cost-effectiveness, it has been used in surveys  \citep{neyman1938}, epidemiological studies  \citep{walker1982, white1982, prentice1986case,  borgan2000exposure, cai2004sample,Kang2009, gray2009weighted}, and   econometrics\citep{manski1977}. As the traditional data collection methods, e.g. questionnaire, is being digitalized and the low-cost sensors become ubiquitous, the construction of a large phase I sample is easier than ever. At the same time, measuring only a small subsample for expensive covariates in early discovery research is likely unavoidable. Two-phase sampling provides a principled approach to integrate large, inexpensive data with small, expensive data. Thus, its adoption in wider science communities will improve the efficiency of  data-based scientific studies, and a general theoretical framework with tools to develop various methods for analyzing two-phase sampling data will be an essential companion for its adoption. 

Unlike two-phase sampling  design that  improves estimation efficiency by thoughtfully selecting samples, auxiliary data, on the other hand, helps researchers gain prediction or inference precision by adding new relevant information into the estimation procedure. It has been used  in  the estimation of  population totals, distribution functions \citep{rao1994estimating}, failure time regression\citep{zhou1995auxiliary}, and so on. The large number of variables we are able to collect today provides us even more opportunities to use this approach. In this paper, we consider using a large amount of auxiliary data collected for the phase I sample in a two-phase sampling data to improve inference and prediction precision based on a possible semiparametric model.  

Take a biomarker epidemiological study as an example.  Epidemiologic cohort studies often follow a large number of participants for many years for the occurrence of disease. If every covariate is measured for every participant, the study can be highly expensive, especially when a disease is rare or the measurement of a covariate, e.g., a new biomarker, requires a new technology.  Two-phase sampling can be implemented to save costs. The expensive covariates are only measured for subjects who subsequently developed a disease and a small cohort random sample while all the other information are measured for the entire cohort. If one is interested in studying association between the biomarker and the disease, a semiparametric survival model such as the Cox regression model is a common choice for random sampling data. The study of association is then translated to a problem of estimating the regression parameters of this model. However, when data are collected with two-phase sampling, how does one estimate the regression parameters?   When auxiliary data such as disease history are collected for all cohort members, how to incorporate this information to improve the estimation precision for the regression parameters?  If our assumed model is wrong, how does one estimate and interpret the estimator under model misspecification? In personalized medicine, biomarkers are considered for individual risk prediction. Risk prediction requires joint estimation of both the parametric and nonparametric parts of the model if a semiparametric model is used. How does one develop such a  prediction model with two-phase sampling data? Solving these problems one after another will be time-consuming. Alternatively, I propose a system to solve these related problems systematically. 

To build such a system, I carefully select and consolidate results from survey sampling \citep{neyman1938m, deville1992calibration}, Z-estimation\citep{ huber1967, pollard1985new}, and the modern empirical process theory\citep{van2000weak}, and then build upon. I make the following three choices: 1) connect a two-phase epidemiological study to survey sampling and adapt general methods from survey methodologies; 2) define a parameter as a functional in $l^{\infty}$ space; 3) use Z-estimation for i.i.d data. It turns out these three choices are suitable for each other. They integrate together and make the framework systematic. Therefore I also call this framework a Z-estimation system.

Several  Z-estimation theorems have been proposed for model development with two-phase sampling. \citet{breslow2007weighted} first used Z-estimation to develop a weight likelihood estimation method for semiparametric models with two-phase sampling. \citet{saegusa2013} provided new Z-theorems for weighted likelihood estimation under sampling without replacement at the Phase II sampling step. \citet{breslow2008} developed Z-theorem with estimated nuisance parameters of possibly infinite dimensional parameter, and \citet{nan2012} further allowed the nuisance parameter to be a function of a Euclidean parameter of interest and convergence of nuisance parameters at a non $\sqrt{n}$-rate. In this paper, I extend the general theory for two-phase sampling proposed by \citet{breslow2007weighted} and expanded it in several directions. I include both likelihood  and non-likelihood based inference to accommodate broader types of estimators such as  moment estimators,  least-squares estimators and so on.  I  allow for model misspecification. I also consider joint estimation of both Euclidean and non-Euclidean for prediction problems. %In addition, I simplify the use of Z-estimation and empirical processes theory results so that the framework could be easily adopted by applied statisticians.  At last, we selected and built tools that can facilitate the reuses of existing results and moved toward a systematic development of a series of related statistics at one time. 

As the first paper on this systematic approach we  focus on developing and demonstrating the system within a limited scope. We select several classical estimation problems that were previously considered separately. We will show a procedure to  develop  them all at once systematically.   The scope of our study includes
\begin{enumerate}
	\item  estimating the same parameter with  a simple random sample,  a more informative subsample generated from a two-phase sampling design(),  and a combined data of this subsample with additional auxiliary information.  %This is a common problem with big data where data  with different underlying generating processes are often assembled together. In practice, this assembly rarely happens once  since new data sources are discovered overtime and then added to the existing ones. 
   	\item   estimating Euclidean, non-Euclidean parameters , and some combinations of these two types of parameters from  the same model; %For example, estimating a regression parameter to study  how a risk factor is associated with a disease from a survival model is considered a related estimation problem to the prediction of an individual's risk of the disease based on the same model.
   	\item  estimating with and without model misspecification;
\end{enumerate}
   as well as any combinations of these estimation problems.

    This scope  should be treated only as a precursor of the estimation problems our proposed system could handle simultaneously. We believe once the framework of this system is established,  the scope could be gradually expanded to allow simultaneously methods development for more data  types and inference problems from similar but different models. Once the system is clearly demonstrated, we encourage the research community to extent the scope of this system.

      	The paper is organized as follows. In section 2, we summarize key concepts repeatedly used in this paper and previous methods proposed for solving individual problems. We  first introduce various data structures this system concerns. Next, we explain  the procedure in this Z-estimation system and a set of theoretical tool for establishing consistency and asymptotic variance.   Then we illustrate  this system by applying it to the additive hazards model. Finally, we summarize the strength and the limitation of this system.

\section{A Z-estimation System}\label{sec:Z-system}
Let us now discuss the proposed estimation system. We start with the parameter this estimation system is concerned with and how we conceptualize it.

\subsection{Parameter  as a functional}\label{sec:para}

In our Z-estimation system we consider the (target) parameter of our estimator as a functional associated with a true underlying distribution  $P$ of the observed data $X$.  This functional is ...  More specifically, let $\mathbb{L}_0$ denote the parameter space. Parameter $\alpha$ is defined as a functional 
\[
\alpha: \mathcal{P}\mapsto \mathbb{L}_0
\]
that is implicitly identified by 
\begin{equation}\label{def:para_equation}
\Psi(\alpha) = 0
\end{equation}
and $\Psi(\alpha)$ can be written as 
\[
\Psi(\alpha)= P\psi_{\alpha}(X) =  \int \psi(X)dP
\] 
Throughout this paper, we use a linear operator symbol to stand for taking expectation.  $\psi_{\alpha}(X)$ can be a score function derived from an assumed model $\{P_{\alpha}, \alpha \in \mathbb{L}_0\}$ or a pseudoscore function motivated from  the assumed model. It can also be motivated from non-likelihood approach such as the method of  least squares.

The true parameter $\alpha_0$ is defined as the value of the  functional evaluated at the true distribution $P$: 
\begin{equation}\label{def:para_functional}
\alpha_0=\alpha(P).
\end{equation}
We require $\alpha_0$ to be unique.  Moreover, we require  both $\alpha$ and  $\psi_{\alpha}(X)$ belong to $ l^{\infty}(H)$ where $H$ is an arbitrary index set.  %	In addition,  we require that both $\psi_{\alpha}(X)$ and $P\psi_{\alpha}(X)$ are uniformly bounded. 

\begin{remark}
	When $\psi_{\alpha}(X) \in l^{\infty}(H)$, \eqref{def:para_equation} is  equivalent to a collection of equations in $\mathbb{R}$:
	\begin{equation}\label{def:para_equation_h}
	\Psi(\alpha)h= P\psi_{\alpha}(X)h= P\psi_{\alpha,h}(X)=0 \mbox{ for every } h \in H. 
	\end{equation} We see from this display that each $h$ indexes an  equation in real line, for which reason $H$ is called an index set. The size of this collection can be infinite depending on the size of $H$. In consequence,  the parameter $\alpha$ defined through $\Psi(\alpha)$ is allowed to be  finite or infinite-dimensional. 
	
	For example, an Euclidean parameter $\theta \in \mathbb{R}^p$ is considered in this system since $\theta$ can also be identified  as a map that maps $h_1$ in the unit ball in $\mathbb{R}^p$ to the real line as $ \theta h_1 = h_1^T\theta$. Some non-Euclidean parameter  is included as well. Consider the cumulative hazard $\Lambda$ over a time interval $[0, \tau]$.  If we treat $\Lambda$ as a map  defined by  $\Lambda h_2 = \int_0 h d\Lambda$, then $\Lambda \in l^{\infty}(H_2)$  where  $H_2$ is the set of bounded functions over $[0,\tau]$.  In addition,  the combination of  Euclidean and non-Euclidean parameters can also be treated as a single parameter and studied by this system.  Let  $H$ = $(H_1, H_2)$. Then ($\theta$, $\Lambda$) is uniquely identified by $\alpha \in l^{\infty}(H)$ where $\alpha h=h_1^T\theta+\int_0^\tau h_2 d\Lambda$.  Since we are   able to consider Euclidean and non-Euclidean parameters together as a single parameter, we can  also estimate them simultaneously in this system.
	
	%	By considering space $\mathbb{L}= l^{\infty}(H)$, on one hand we concisely describe an infinite number of equations with one symbol $\psi_{\alpha}$, on the other hand the index $h$ facilitates us to focus on a single dimensional function $\psi_{\alpha,h}$ at one time in this infinite-dimensional problem. In later sections we will see an abundance of convenience brought by restricting $\mathbb{L} = l^{\infty}(H)$. 
\end{remark}

\begin{remark}
	%	Model parameter was often defined as a component of  an assumed model $\{P_{\alpha}, \alpha \in \mathbb{L}_0\}$. Unlike this type of definition, our parameter  is defined by associating it with a non-parametric model $\mathcal{P}$.  We do not have restriction on $\mathcal{P}$ except for the regularity condition that $PXX^T <\infty$ for $P \in \mathcal{P}$. Our definition is inspired by  \citet{newey1994asymptotic} and \citet{binder1992}.
	
	A model parameter is often defined as a component of  an assumed model $\{P_{\alpha}, \alpha \in \mathbb{L}_0\}$. By allowing the true distribution $P$ for $X$ to vary in $\mathcal{P}$, we define the target of our estimators allowing for model misspecification.  When $P \notin \{P_{\alpha}, \alpha \in \mathbb{L}_0\}$,  we still have an interpretation of the parameter according to our definition and this  parameter will be consistently estimated as shown later under some conditions.  When $P \in \{P_{\alpha}, \alpha \in \mathbb{L}_0\}$, the parameter from our definition agrees with the model parameter in $P_{\alpha}$.  Our definition is inspired by  \citet{newey1994asymptotic} and \citet{binder1992}.
	
	%An assumed model $\{P_{\alpha}, \alpha \in \mathbb{L}_0\}$ can motivate  $\psi_{\alpha}$ but do not determine  $\psi_{\alpha}$ completely.
	
\end{remark}

%When we set the parameter space to be $\mathbb{L}_0 \subset l^{\infty}(H)$ and  $H$ to be arbitrary, we are concerned with a large range of parameters.  Any parameter that can be uniquely identified by a map $\alpha$ in $l^{\infty}(H)$  is considered in this system. 
%This includes Euclidean parameter $\theta \in \mathbb{R}^p$ since $\theta$ can be uniquely identified  by the map $\alpha$ in $l^{\infty}(H_1)$   where $\alpha h_1 = h_1^T\theta$ and  $H_1$ is the unit ball in $\mathbb{R}^p$. Some non-Euclidean parameter  such as the cumulative hazard $\Lambda$ defined on a time interval $[0, \tau]$ is also considered, because $\Lambda$ can be uniquely identified  by the map $\alpha$ in $l^{\infty}(H_2)$  where $\alpha h_2 = \int_0 h d\Lambda$ and $H_2$ is the set of bounded functions over $[0,\tau]$.  We can also treat parameters of these two  types such as ($\theta$, $\Lambda$) as a single parameter in this system. Let  $H$ be the product space of $H_1$ and $H_2$. ($\theta$, $\Lambda$) is then uniquely identified by $\alpha \in l^{\infty}(H)$ where $\alpha h=h_1^T\theta+\int_0^\tau h_2 d\Lambda$.  These three examples demonstrate the large range of parameters we are concerned with in this Z-estimation system. The fact that we are able to consider Euclidean and non-Euclidean parameters together as a single parameter provides a possibility and convenience to jointly estimate them simultaneously. 

\section{Data Structures and Assumptions}\label{sec:sampling}
Suppose we are interested in estimating a parameter $\alpha$ that describes an aspect of a population. Let $X$  denote the characteristics of an individual useful for estimating $\alpha$. 

\textbf{With random sampling (RS)}, the follow assumption is made about the data $X_1, X_2, \dots, X_N$:
	\begin{assumption}\label{assum: i.i.d.}
	$N$ individuals are randomly selected from a population and vectors $X_1,\dots,X_N$ are i.i.d. random samples from a probability distribution $P$ on a measurable space $(\mathcal{X},\mathcal{A})$. $P \in \mathcal{P}$, which is the set of all possible distributions $P$ for $X$ such that $E(X^TX) < \infty$. 
\end{assumption}

\textbf{With two-phase sampling}, the data is collected in two steps. The vector of interest $X$ is also divided into two parts $X^{I}$ and $X^{II}$ in which $X^{II}$ is often costly.    At phase I we  randomly select $N$ individuals from a population  and only measure  $X^{I}$.  Moreover, we may collect some auxiliary variables $U$ that are available for all $N$ individuals. 

Next, based on $V= (X^I,U) \in \mathcal{V}$, we independently generate a phase II selection indicator $R$ for each individual from a Bernoulli distribution with the probability of $P(R=1|V) =\pi_0(V)$, in which $\pi_0$ is an a priori function of $V$. Based on $R$ we construct our phase II subsamples and only these individuals are sent for the measurement of   $X^{II}$.

We call this sampling design two-phase variable probability sampling (VPS). When $V$ indicates strata and $\pi_0(V)$ is a constant function for each stratum, this sampling method becomes the familiar stratified sampling and VPS 1 described in \citet{kalbfleisch1988}. However, we consider a more general two-phase sampling design that allows $V$ and $\pi_0(V)$ to be continuous. 
  
 	Two-phase VPS generates  data in the form of $(R_i,V_i, X_i^{II}R_i), i=1, \dots,N$.  For this data, we make the following assumptions in addition to  A.\ref{assum: i.i.d.}:	
	\begin{assumption}\label{assum: MAR}
			Whether or not a subject is selected into phase II depends only on the observed data $V$ at phase I: 
			$P(R_i=1|X_i,U_i) =P(R_i=1|V_i)= \pi_0(V_i)$, in other words, missing at random(MAR) \citep{rubin1976inference}.
	\end{assumption}
	
	\begin{assumption}\label{assum: bounded_below}
	\mbox{Sampling probability is bounded below from 0:}
	$\pi_0(v) \geq \sigma > 0$ for all $v \in \mathcal{V}$.
	\end{assumption}

	\begin{assumption}\label{assum: i.i.d.-2ph}
	$(R_i, V_i, X_i^{II}R_i)$ are i.i.d. following a probability distribution $Q$ on a measurable space $(\mathcal{R}\times \mathcal{V} \times \mathcal{X}, \mathcal{B})$. $Q \in \mathcal{Q}$, which is the set of all possible distributions for $(R, V, X)$ such that $E[(R, V, X)^T(R, V, X)] < \infty$.
	\end{assumption}
	
	Note by two-phase VPS, $R_i$ is independently generated. If we assume $V_i, i =1, 2, \dots, N$ are i.i.d., under the assumption of A.\ref{assum: i.i.d.} - \ref{assum: MAR}, $(R_i, V_i, X_i^{II}R_i)$ are ensured to be i.i.d.. A.\ref{assum: MAR} - \ref{assum: bounded_below}  are also assumed in \citet{robins1995} for their missing data problems. The assumptions of $P \in \mathcal{P}$ and $Q \in \mathcal{Q}$ in A. \ref{assum: i.i.d.} and \ref{assum: i.i.d.-2ph} are similar to \citet{newey1994asymptotic} and allow us to estimate under general misspecification.

So far $V$ is only used for deciding the phase II membership $R$. These information may also be incorporated for precision improvement  at the estimation stage. Let vector $\tilde{V}=\tilde{V}(V)$ of q-dimension be the phase I quantity we choose to incorporate. It can be a part of $V$  or a function of the variables in $V$. 

\textbf{When incorporating auxiliary variables $\tilde{V}$} to further  improve  estimation precision,  in addition to A.\ref{assum: i.i.d.} - A.\ref{assum:  i.i.d.-2ph}, we assume:
	\begin{assumption}\label{assum: vtilde-bounded}
	$\tilde{V}$ is bounded.
	\end{assumption}

	\begin{assumption}\label{assum: vtilde-posdef}
	$Q\tilde{V}\tilde{V}^T$ is positive definite. 
	\end{assumption}

\subsection{Z-Estimators Construction}\label{sec:estimator}
To develop  various estimators  systematically,  every estimator we propose is a  Z-estimator regardless of the data types, which is  defined as a solution to some estimating equations (EE).  To simplify and systematize the asymptotic analysis of these estimators, we are only concerned with  those Z-estimators from the EE of i.i.d. functions. 
%We consider the ordinary EE   and  inverse probability weighted EE (IPW-EE)  of i.i.d. functions for problems 1 and 2 respectively. In problem 3,  we use  a calibration technique \citep{deville1992calibration}  to adjust the weights used  in the IPW-EE  based on some phase I auxiliary information. It turns out the problem of weight calibration  can be transformed to a problem of solving additional EE of i.i.d. functions of these auxiliary information. The solution to problem 3 is still a Z-estimator.

Let $\mathbb{P}_N$ stand for taking expectation  under the empirical measure $\mathbb{P}_N:
\mathbb{P}_N f= \frac{1}{N} \sum_{i=1}^N f(X_i)$.  \textbf {When data on $X$ are obtained from RS}, 
%we consider a  map $\Psi_N : \mathbb{L}_0 \mapsto \mathbb{L}$ of the form $\Psi_N (\alpha) = \mathbb{P}_N\psi_{\alpha}(X)$ where
%$\mathbb{L} = l^{\infty}(H)$. Our construction  of EE is as follows.  
%	\begin{equation}\label{EE:complete}
%\Psi_N (\alpha)=0
%\end{equation}
%or equivalently 
we propose our RS estimator $\hat{\alpha}$ as the solution to the EE:
\begin{equation}\label{EE:complete}
\mathbb{P}_N\psi_{h, \alpha}(X) = 0, \quad \forall  h \in H.
\end{equation}

\textbf{When data on $X$ are obtained from two-phase VPS}, $X$ is complete only for the phase II subsample. Let $\psi_{\alpha}^*(X,V,R)=\frac{R}{\pi_0(V)}\psi_{\alpha}(X)$. 
% We consider a new random map $\Psi_N^*: \mathbb{L}_0 \mapsto \mathbb{L}$ given by $\Psi_N^* (\alpha) =\mathbb{P}_N\psi_{\alpha}^*(X,V,R)$. 
We propose our two-phase VPS estimator $\hat{\alpha}^* $ as the solution to an inverse probability weighted EE (IPW-EE):
\begin{equation}\label{EE:twophase}
%\Psi_N^*(\alpha)=0.
\mathbb{P}_N\psi_{\alpha}^*(X,V,R)=\mathbb{P}_N\frac{R}{\pi_0(V)}\psi_{\alpha}(X) = 0, \quad \forall   h \in H. 
\end{equation}
%Recall R is the indicator of whether a phase I member is selected into the phase II subsample,  only the member of which contributes to the construction of $\hat{\alpha}^*$. 
%because the contribution from each subject to the overall estimating equation is weighted by the inverse of this subject's selection probability. 

\textbf{When data on $X$ from phase II and data on $\tilde{V}$ from phase I in two-phase VPS are combined}, we propose a new estimator that could potentially improve the estimation precision compared to using phase II data alone. We modify the weights in \eqref{EE:twophase} for $\hat{\alpha}^*$ so that the totals of $\tilde{V}$ calculated based on the phase I observations are exactly estimated by their phase II estimates with some modified weights\citep{deville1992calibration}.  Let
\begin{align*}
\psi^{c}_{1,\alpha,\gamma}(X,V,R) &= \frac{R}{\pi_0(V)}exp(-\gamma^T\tilde{V})\psi_{\alpha}(X), \\
\psi^{c}_{2,\gamma}(V,R) & = \frac{R}{\pi_0(V)}exp(-\gamma^T \tilde{V}) \tilde{V} -\tilde{V}
\end{align*}
where $\gamma \in \Gamma$ and $\Gamma \subset \mathbb{R}^q$. %We consider a new random map $ \Psi^{c}_N : \mathbb{L}_0 \times \Gamma \mapsto \mathbb{L}\times \mathbb{R}^q$ given by $\Psi^{c}_N (\alpha, \gamma) = (\Psi^{c}_{N,1}(\alpha, \gamma), \Psi^{c}_{N,2}(\alpha, \gamma)) = (\mathbb{P}_N\psi^{c}_{1,\alpha,\gamma}(X,V,R), \mathbb{P}_N\psi^{c}_{2,\alpha,\gamma}(V,R))$.  
We propose our  calibrated two-phase VPS estimator  $(\hat{\alpha}^{c},\hat{\gamma})$ as the solution to the following EE:
\begin{equation}\label{EE:calibrated_estimators}
%\Psi^{c}_N(\alpha,\gamma)=0.
\begin{split}
\mathbb{P}_N\psi^{c}_{1,\alpha,\gamma,h}(X,V,R) =  \mathbb{P}_N\frac{R}{\pi_0(V)}exp(-\gamma^T\tilde{V})\psi_{ \alpha,h}(X) =0, \quad \forall h \in H\\
\mathbb{P}_N\psi^{c}_{2,\alpha,\gamma}(V,R) =  \mathbb{P}_N \left\{\frac{R}{\pi_0(V)}exp(-\gamma^T \tilde{V}) \tilde{V} -\tilde{V}\right\} = 0.
\end{split}
\end{equation}
Comparing to EE in \eqref{EE:twophase}, the additional information on $\tilde{V}$ from phase I is incorporated into the new estimator through the adjusted weight $\frac{1}{\pi_0(V)}exp(-\gamma^T\tilde{V})$, and this additional information may  improve our estimation efficiency.

The above  EE are  motivated as follows. Let Poisson deviance
$G\{w, \frac{1}{\pi_0(v)}\} = wlog\{\frac{w}{1/\pi_0(v)}\}-\{w -\frac{1}{\pi_0(v)}\}$ as our distance measure for comparing the adjusted weight $w$ to the original weight $1/\pi_0(v)$. According to the calibration criteria,  our goal is to  find a new series of weights $w_i, i=1,2, \cdots, N$ such that $\sum_{i=1}^N R_i G\{w_i,\frac{1}{\pi_0(v_i)}\}$ is minimized and at the same time  these weights are subject to the constraints
\begin{equation}\label{constraints}
\sum_{i=1}^{N}\tilde{v_i} =\sum_{i=1}^N R_i w_i \tilde{v_i}.
\end{equation}
We use the Lagrange multipliers method. Let $\gamma \in \Gamma $ be a q-vector of Lagrange multipliers corresponding to the constraints \eqref{constraints} and let $f(w_i, \gamma)$ be the Lagrange functions:
\[
f(w_i, \gamma) = G\{w_i, \frac{1}{\pi_0(v_i)}\} + \gamma^T( \sum_{i=1}^{N}\tilde{v_i}-\sum_{i=1}^N R_i w_i \tilde{v_i}) 
\]
for $i =1,2,\dots,N$. Solving 
\[
\nabla _{w_i}f(w_i, \gamma) =0 \mbox{ and } \nabla _{\gamma}f(w_i,\gamma) =0
\]
yields
\begin{equation}\label{new_weight}
w_i = \frac{exp(-\gamma^T \tilde{v_i})}{\pi_0(v_i)}.
\end{equation} 
Replacing the original weight $1/\pi_0({v_i})$ in \eqref{EE:twophase} by $w_i$  in \eqref{new_weight} yields the  EE  for $\hat{\alpha}^c$: $ \mathbb{P}_N\psi^{c}_{1,\alpha,\gamma}(X,V,R) = 0$. Combining  \eqref{constraints} and \eqref{new_weight} yields the EE for $\hat{\gamma}$: 
\begin{equation}\label{new_ee}
\sum_{i=1}^{N}\tilde{v_i} = \sum_{i=1}^{N} \frac{R_i}{\pi_0(v_i)} exp(-\gamma^T \tilde{v_i}) \tilde{v_i}.
\end{equation}

%Our development of  $(\hat{\alpha}^c, \hat{\gamma})$ also illustrates why we are interested in expanding from  likelihood-based inference in \citet{breslow2007weighted} to non-likelihood based inference. In this case, the estimator for $\gamma$ is motivated from the Lagrange multipliers method rather than a likelihood function.

\begin{remark}
An alternative approach to incorporate phase I information is to estimate the phase II selection probability by a prediction model using phase I covariates. We prefer the calibration method because it fits right into the Z-estimation system. Since calibrated weights can be written as a solution to  EE of i.i.d. functions as shown in \eqref{new_ee} and  EE is  allowed to be  infinite dimensional in this system, adding a new set of EE for finding the calibrated weights to the existing EE does not introduce extra complexity. In addition, unlike the prediction model approach, the calibration method does not dependent on an assumed model. It is compatible to the Z-estimation system that is model assisted but not model dependent. 
\end{remark}

%We consider $\gamma$ ranges over the parameter space $\Gamma$ where $\Gamma \subset \mathbb{R}^q$ and $(\alpha, \gamma)$ ranges over the parameter space $ \mathbb{L}_0 \times \Gamma$. %Recall %$\mathbb{L}_0$ is the parameter space for $\alpha$,$\mathbb{L}_0 \subset \mathbb{L}$ and $\mathbb{L}=l^{\infty}(\mathcal{\part{title}H})$. %Let $\mathbb{T} = \mathbb{L}\times \mathbb{R}^q$. %Then $\mathbb{T}_0 \subset \mathbb{T}$ 

\subsection{Conditions as Components}\label{sec:conditions}

$\hat{\alpha}$, $\hat{\alpha}^{*}$ and $\hat{\alpha}^{c}$ all come from zero-valued averages of  i.i.d. functions. These functions  are  closely related through the common map $\psi_{\alpha}(X)$. Hence  we consider transforming the asymptotic studies of these three estimators to a sequence of analyses on $\psi_{\alpha}(X)$. Below is a list of conditions on $\psi_{\alpha}$ that is necessary for proving the asymptotic normality of these three estimators at a $\sqrt{N}$ rate.

\begin{condition}
	For each $\alpha \in \mathbb{L}_0 $ of a normed space and every $h$ in an arbitrary set $H$, let $x \mapsto \psi_{\alpha,h}(x)$ be a measurable function, such that 
	\begin{itemize}
		\item[\circled{1}]  $\mathcal{G}\equiv\{\psi_{\alpha,h}, \alpha \in \mathbb{L}_0, h \in H\}$ is Glivenko-Cantelli.
		\item[\circled{2}] $P \psi_{\alpha}(\alpha_0)=0$ and for every $\epsilon > 0,$
		\[
		\underset{\alpha:\|\alpha-\alpha_0\|\geq \epsilon}{\inf}\|P \psi_{\alpha}-P \psi_{\alpha_0}\|_H >0.
		\] 
		\item[\circled{3}] The class $\mathcal{F} \equiv \{ \psi_{\alpha,h}: \|\alpha-\alpha_0\|< \delta, h \in H\}$, with finite and integrable envelope function, is $P$-Donsker for some $\delta>0$.
		\item[\circled{4}] As a map into $l^\infty (H),$ the map $\alpha \mapsto P\psi_{\alpha}$ is Fr\'echet-differentiable at $\alpha_0$, with a derivative $\dot{\Psi}_0: lin \mathbb{L}_0 \mapsto l^\infty ( H)$ that has a continuous inverse on its range.
		\item[\circled{5}] $\|P(\psi_{\alpha} - \psi_{\alpha_0})^2 \|_{H} \rightarrow 0$ as $\alpha \rightarrow \alpha_0$.
	\end{itemize}
\end{condition}
If $\alpha$ is finite dimensional, the Euclidean norm is used.  If it is infinite dimensional,  the super norm  is considered. Denote by $\|\cdot\|_H$ the sup norm $\|x\|_H= sup_{h \in H}|x(h)|$. For a product space, the norm is defined as a sum $\|(x,y)\|=\|x\|+\|y\|$.
Note since  in this system $\psi_{\alpha}$ does not need to be a score function as in \citep{breslow2007weighted} from a likelihood function and it can just be motivated from a model, some technical challenges in verifying the above conditions of $\psi_{\alpha}$ can be avoided by choosing $\psi_{\alpha}$ wisely.
%{\color{red}{Why in thm 19.26 of \citep{vdv1998}, condition 3 only requires a finite envelope function rather than both finite and integrable?}}.
\subsection{Theoretical Tools and Asymptotic Results}\label{sec:tools}
In this section, we develop a series of tools that can be used for quickly establishing the consistency and the limiting distributions of the proposed estimators  once the analytic conditions in the proceeding section are established.  We start with a  consistency theorem for very general Z-estimators. %and then three consistency corollaries for our Z-estimators $\hat{\alpha}$, $\hat{\alpha}^{*}$ and $\hat{\alpha}^{c}$. 

\begin{theorem}\label{thm:consistency}
	Let $\Psi_N: \mathbb{L}_0 \mapsto \mathbb{L}$ be a random map and $\Psi: \mathbb{L}_0\mapsto \mathbb{L}$ be a deterministic map, both from the parameter space $\mathbb{L}_0$, a subset of Banach space, to the  Banach space $\mathbb{L}$. Let $||\cdot||_{\mathbb{L}_0}$ and  $|| \cdot||_{\mathbb{L}}$ be the norms for $\mathbb{L}_0$ and $\mathbb{L}$ respectively. 
	Consider
	\begin{conditionA}
		\[
		\underset{\alpha \in \mathbb{L}_0}{\sup}\|\Psi_N(\alpha)-\Psi(\alpha)\|_\mathbb{L} \overset{p}\rightarrow 0;\]
	\end{conditionA}
	\begin{conditionB}
		$\Psi(\alpha_0)=0$ and for every $\epsilon > 0,$
		\[
		\underset{\alpha:\|\alpha-\alpha_0\|\geq \epsilon}{\inf}\|\Psi(\alpha)-\Psi(\alpha_0)\|_\mathbb{L} > 0.
		\] 
	\end{conditionB}
	If both conditions are satisfied, then $\hat{\alpha}$ that satisfies $\Psi_N(\hat{\alpha})=o_p(1)$ converges to $\alpha_0$ in probability.
\end{theorem}

\begin{proof}
	Based on condition A,
	\begin{equation}\label{op1_1}
	\|\Psi(\hat{\alpha})-\Psi_N(\hat{\alpha})\|_\mathbb{L} \leq \underset{\alpha \in \mathbb{L}_0}{\sup} \| \Psi(\alpha)-\Psi_N(\alpha)\|_\mathbb{L}= o_p(1). 
	\end{equation}
	By assumption, we have $\Psi(\alpha_0)=0$ and $\Psi_N(\hat{\alpha}) =o_p(1)$, so
	\begin{equation}\label{op1_2}
	\|\Psi_N(\hat{\alpha})-\Psi(\alpha_0)\|_{\mathbb{L}} = o_p(1).
	\end{equation}
	Combining \eqref{op1_1} and \eqref{op1_2} results in
	\[
	\|\Psi(\hat{\alpha})-\Psi(\alpha_0)\|_{\mathbb{L}}\leq \|\Psi(\hat{\alpha})-\Psi_N(\hat{\alpha})\|_{\mathbb{L}}+\|\Psi_N(\hat{\alpha})-\Psi(\alpha_0)\|_{\mathbb{L}}
	=o_p(1),\]
	which means for any $\eta >0$,
	$P(\|\Psi(\hat{\alpha})-\Psi(\alpha_0)\|_{\mathbb{L}} \geq \eta)\rightarrow 0$ as $N \rightarrow \infty$. 
	
	According to condition B, for any $\epsilon > 0$, there exists $\eta >0$ such that event $A \equiv \{
	\|\hat{\alpha}-\alpha_0\|_{\mathbb{L}_0} \geq \epsilon >0\}$ implies event $B \equiv \{ \|\Psi(\hat{\alpha})-\Psi(\alpha_0)\|_{\mathbb{L}} \geq \eta\}$, in other words, $A \subset B$, so
	\[
	P(\|\hat{\alpha}-\alpha_0\|_{\mathbb{L}_0}\geq \epsilon ) 
	\leq P( \|\Psi(\hat{\alpha})-\Psi(\alpha_0)\|_\mathbb{L} \geq \eta).
	\]
	Therefore $\forall \epsilon >0$, $P(\|\hat{\alpha}-\alpha_0\|_{\mathbb{L}_0} \geq \epsilon ) \rightarrow 0$ as $N \rightarrow \infty$. We conclude $\|\hat{\alpha}-\alpha_0\|_{\mathbb{L}_0} \rightarrow_p 0$.
\end{proof}

This theorem considers maps between general Banach spaces $\mathbb{L}_0$ and $\mathbb{L}$. $\mathbb{L}$ is not restricted to $l^{\infty}(H)$.  Condition A is often challenging to verify and there is not a standard approach to prove it. In the case of i.i.d. observations with  $\mathbb{L} = l^{\infty}(H)$ as we require for our Z-estimation system, however,   condition A is simplified. Because
\begin{equation}\label{eqn:GC}
\underset{\alpha \in \mathbb{L}_0}{\sup}\|\Psi_N(\alpha)-\Psi(\alpha)\|_\mathbb{L}= \underset{\alpha \in \mathbb{L}_0}{\sup}\|\mathbb{P}_N\psi_{\alpha}-{P}\psi_{\alpha}\|_{H}
=\underset{\alpha \in \mathbb{L}_0, h \in H}{\sup}|\mathbb{P}_N\psi_{\alpha,h}-{P}\psi_{\alpha,h}|,
\end{equation}
we only need to verify condition  \circled{1} in section \ref{sec:conditions}  that $\mathcal{G}\equiv\{\psi_{\alpha,h}, \alpha \in \mathbb{L}_0, h \in H\}$  is Glivenko-Cantelli (GC). As a result, many existing tools about GC classes in  the empirical processes theory literatures  can be immediately used for establishing consistency in our Z-estimation system as shown in the following three consistency corollaries for  $\hat{\alpha}$, $\hat{\alpha}^*$ and $\hat{\alpha}^c$.

% Laster we will see GC condition is effortlessly carried to the new two-phase VPS scenario, the establishment of consistency for $\hat{\alpha}^*$ and $\hat{\alpha}^c$ are simplified together with $\hat{\alpha}$, as shown below. 	The proofs demonstrate  how existing mathematical results on GC classes  extends the consistency property of a simple statistic to a more complex one effort

\begin{corollary} (Consistency of RS estimators)\label{thm: consist-srs}
	If $\psi_{\alpha}$ is a map satisfying conditions \circled{1} and \circled{2}, then
	\[
	\hat{\alpha} \rightarrow_p \alpha_0.
	\] 
\end{corollary}
\begin{proof}
	Let $ \Psi_N(\alpha) =\mathbb{P}_N\psi_{\alpha}$ and $\Psi(\alpha) = {P}\psi_{\alpha}$.
	Under condition \circled{1}, \\ $\underset{\alpha \in \mathbb{L}_0, h \in H}{\sup}|\mathbb{P}_N\psi_{\alpha,h}-{P}\psi_{\alpha,h}|\rightarrow_p 0$, and according to \eqref{eqn:GC},
	$\underset{\alpha \in \mathbb{L}_0}{\sup}\|\Psi_N(\alpha)-\Psi(\alpha)\|_\mathbb{L} \rightarrow_p 0$. Condition A in theorem   \ref{thm:consistency} is therefore satisfied. Under   condition \circled{2}, Condition B is automatically  satisfied. 	By the definition of $\hat{\alpha}$ in \eqref{EE:complete}, $\Psi_N(\hat{\alpha})=0$.  It follows from theorem \ref{thm:consistency} that $\hat{\alpha}\rightarrow_p \alpha_0$. 
\end{proof}

\begin{corollary}(Consistency of two-phase VPS estimators)\label{thm: consist-2ph}\\
	Suppose A \ref{assum: MAR}  - \ref{assum: bounded_below} hold. If $\psi_\alpha$ is a map  satisfying  conditions \circled{1} and \circled{2} with and $\psi_\alpha$ is uniformly bounded,  then 
	\[
	\hat{\alpha}^* \rightarrow_p \alpha_0.
	\] 
\end{corollary}

\begin{proof} We show consistency of $\hat{\alpha}^*$ by adapting theorem \ref{thm:consistency} to $\hat{\alpha}^*$. Let  $ \Psi_N(\alpha) =\mathbb{P}_N\psi^*_{\alpha} = \mathbb{P}_N\frac{R}{\pi_0(V)}\psi_{\alpha}$ and $\Psi(\alpha) = {Q}\psi^*_{\alpha}$. Since $\pi_0(V)$ is bounded from 0 by assumption A. \ref{assum: bounded_below}, $\{\frac{R}{\pi_0(V)}\}$ is a finite class of a square integrable function.  Hence, $\{\frac{R}{\pi_0(V)}\}$ is GC \citep{van2000weak} and it is bounded. Given condition \circled{1}  and our assumption,    $\mathcal{G}$ is also a GC class and  it is uniformly bounded.  Applying the preservation theorem of GC \citep{van2000},  we obtain $\mathcal{G}' = \{\psi^*_{\alpha,h}(R,V,X) = \frac{R}{\pi_0(V)}\psi_{\alpha, h}(X): \alpha \in \mathbb{L}_0, h\in H\}$ is GC. Therefore, condition A in theorem \ref{thm:consistency} is satisfied when $\Psi_N(\alpha) =\mathbb{P}_N\psi^*_{\alpha}$ . Next,  we verify condition B in  theorem \ref{thm:consistency} also holds. By A.\ref{assum: MAR} we have
	\begin{equation}\label{eqn:connection}
	Q\psi^*_{\alpha} = E_Q\left[E_Q\left[\frac{R}{\pi_0(V)}|X, U \right]\psi_{\alpha}(X)\right]=P\psi_{\alpha}(X) 
	\end{equation}
   and condition \circled{2} leads to condition B  immediately when $\Psi(\alpha) = {Q}\psi^*_{\alpha}$.
	Finally, by the definition of $\hat{\alpha}^*$ in \eqref{EE:twophase}, $\Psi_N(\hat{\alpha}^*) =0$. It follows from theorem \ref{thm:consistency} that $\hat{\alpha}^* \rightarrow_p \alpha_0$.
	%which claims a class of function formed by the product of two GC classes with integrable envelope functions is still GC. Using this theorem we can have a systematic approach to show GC by showing each component of the product is GC. 
\end{proof}

Before establishing the consistency corollary for $\hat{\alpha}^c$,  we first give two lemmas. Their proofs are provided in Appendix A.
\begin{lemma}\label{lemma:unique_solution}
	Assume $\Gamma$ is a compact convex subset of $\mathbb{R}^q$ with $0$ as an interior point. Assume that A.\ref{assum: MAR} and A.\ref{assum: vtilde-bounded} hold.  If $\alpha_0$ is the unique solution to the equation $P\psi_{\alpha} =0$, then $(\alpha_0, 0)$ is the unique solution to the equations  $Q\psi^{c}_{1,\alpha,\gamma}(X,V,R) =0$  and 
	$Q\psi^{c}_{2,\gamma}(V,R) =0$.
\end{lemma}

\begin{lemma}\label{fgamma}
	Suppose $\Gamma$ is a compact convex subset of $\mathbb{R}^q$ with $0$ as an interior point. Assume that A.\ref{assum: bounded_below} and A.\ref{assum: vtilde-bounded} hold. Let $f_{\gamma}(V,R) = \frac{R}{\pi_0(V)}exp(-\gamma^T\tilde{V})$. Then $\{f_{\gamma}, \gamma \in \Gamma\}$ is uniformly Lipschitz and uniformly bounded. 
\end{lemma}

\begin{corollary} (Consistency of calibrated two-phase VPS estimators)\label{thm: consist-cali}
	Suppose $ \Gamma$ is a compact convex subset of $\mathbb{R}^q $ with 0 as an interior point. Assume that A.\ref{assum: MAR}-\ref{assum: vtilde-bounded} hold. If $\psi_\alpha$ is a map satisfying conditions \circled{1} and \circled{2} and it is uniformly bounded, then
	\[
	\hat{\alpha}^{c}\rightarrow_p \alpha_0.
	\]
\end{corollary}

\begin{proof}
	We prove this corollary by applying theorem \ref{thm:consistency} to $(\hat{\alpha}^{c}, \hat{\gamma})$.  To avoid confusion, we replace notation $\alpha$ and $\alpha_0$  in theorem \ref{thm:consistency} by $\tau$ and $\tau_0$ where  $\tau = (\alpha, \gamma)$ and $\tau_0= (\alpha_0, 0)$. We let $\Psi_N(\tau) =  (\mathbb{P}_N \psi^{c}_{1,\alpha,\gamma},\mathbb{P}_N  \psi^{c}_{2,\gamma})$ and $\Psi(\tau) =  (Q \psi^{c}_{1,\alpha,\gamma}, Q \psi^{c}_{2,\gamma})$. They are the maps from the parameter space $\mathbb{L}_0 \times \Gamma$ to  $\mathbb{L}  \times \mathbb{R}^q$.
	% Now $\mathbb{L} =  l^{\infty}(H)\times \Gamma$ 
	We will prove conditions in theorem  \ref{thm:consistency} still hold after this adaptation.
	
	First,  by lemma \ref{fgamma} $\{f_{\gamma}, \gamma \in \Gamma  \}$ is uniformly Lipschitz and uniformly bounded. It  follows from theorems 2.4.1 and 2.7.11  \citep{van2000weak} that $\{f_{\gamma}, \gamma \in \Gamma \}$ is a  GC class.  By condition \circled{1} and our assumption, $\mathcal{G}$ is also a GC  class and it is uniformly bounded. Applying the preservation theorem for GC classes \citep{van2000} we obtain $\{\psi^{c}_{1,\alpha,\gamma, h}=f_{\gamma}\psi_{\alpha,h},(\alpha, \gamma) \in \mathbb{L}_0 \times \Gamma, h \in H\}$ is GC. By assumption A. \ref{assum: vtilde-bounded},  $\{\tilde{V}\}$ is a finite class of a square integrable function, thus it is GC \citep{van2000}. Applying the preservation theorem again we obtain $\{\psi^{c}_{2,\gamma} = f_{\gamma}(V,R)\tilde{V}-\tilde{V}, \gamma \in  \Gamma\}$ is  GC.  As a result,
	%\begin{align*}
	$\underset{(\alpha, \gamma) \in \mathbb{L}_0 \times \Gamma, h \in H}{\sup}|\mathbb{P}_N \psi_{1,\alpha,\gamma,h}^{c}-Q\psi_{1,\alpha,\gamma,h}^{c} | \rightarrow_p 0$. and 
	%\mbox{ }
	$\underset{\gamma \in \Gamma}{\sup}\|\mathbb{P}_N \psi_{2,\gamma}^{c}-Q \psi_{2,\gamma}^{c}\|_\mathbb{E} \rightarrow_p 0.$ Because
	\begin{multline*}
	\underset{\tau \in \mathbb{L}_0 \times \Gamma}{\sup}\|\Psi_N(\tau)-\Psi(\tau) \|_{\mathbb{L}\times \mathbb{R}^q} \leq \\
	\underset{(\alpha, \gamma) \in \mathbb{L}_0 \times \Gamma}{\sup} \|\mathbb{P}_N\psi_{1,\alpha,\gamma}^{c}-Q\psi_{1,\alpha,\gamma}^{c}\|_{H}  +
	\underset{\gamma \in \Gamma}{\sup}
	\|\mathbb{P}_N \psi_{2,\gamma}^{c}-Q \psi_{2,\gamma}^{c}\|_\mathbb{E},
	\end{multline*}
	we have	$\underset{\tau \in \mathbb{L}_0 \times \Gamma}{\sup}\|\Psi_N(\tau)-\Psi(\tau) \|_{\mathbb{L}\times \mathbb{R}^q} \rightarrow_p 0$. 
	% \begin{align*}
	% &\underset{(\alpha, \gamma) \in \mathbb{L}_0 \times \Gamma}{\sup}\|\Psi^{c}_N(\alpha, \gamma)-\Psi^{c}(\alpha, \gamma) \|_{\mathbb{L}\times \mathbb{R}^q}\\
	% = & \underset{(\alpha, \gamma) \in \mathbb{L}_0 \times \Gamma}{\sup}\left\{\|\Psi^{c}_{N,1}(\alpha, \gamma)-\Psi^{c}_1(\alpha, \gamma) \|_H+\|\Psi^{c}_{N,2}(\alpha, \gamma)-\Psi_2^{c}(\alpha, \gamma) \|_\mathbb{E}\right\}\\
	% \leq &\underset{(\alpha, \gamma) \in \mathbb{L}_0 \times \Gamma}{\sup}\|\Psi^{c}_{N,1}(\alpha, \gamma)-\Psi^{c}_1(\alpha, \gamma) \|_H+\underset{ \gamma \in \Gamma}{\sup}\|\Psi^{c}_{N,2}(\alpha, \gamma)-\Psi_2^{c}(\alpha, \gamma) \|_\mathbb{E}\\
	% =&\underset{(\alpha, \gamma) \in \mathbb{L}_0 \times \Gamma, h \in H}{\sup}|\mathbb{P}_N \psi_{1,\alpha,\gamma}^{c}-Q\psi_{1,\alpha,\gamma}^{c} | + \underset{\gamma \in \Gamma}{\sup}\|\mathbb{P}_N \psi_{2,\gamma}^{c}-Q \psi_{2,\gamma}^{c} \|_\mathbb{E}
	% \rightarrow_p 0.
	% \end{align*}
	
	Next, we show condition B holds for  $\Psi$. By lemma \ref{lemma:unique_solution}, $\tau_0$ is the unique solution to $\Psi(\tau)=0$. We only need to show $\tau_0$ is also a well-separated solution. 
	We divide the set %$\{\tau:\|\tau-\tau_0\|\geq \epsilon\} = 
	$\{\alpha,\gamma:\|\alpha-\alpha_0\|+\|\gamma-0\|\geq \epsilon\}$ into two parts $ A\equiv \{\gamma=0\} \cap \{\alpha,\gamma:\|\alpha-\alpha_0\|+\|\gamma-0\|\geq \epsilon\}$  and $
	B \equiv \{\gamma \neq 0\} \cap \{\alpha,\gamma:\|\alpha-\alpha_0\|+\|\gamma-0\|\geq \epsilon\}$.
	We examine the infima of $\|\Psi(\tau)-\Psi(\tau_0) \|_{\mathbb{L}\times \mathbb{R}^q}$ over sets A and B separately. If both infima are greater than 0, then the infima over $A\cup B$ is great than 0.   By \eqref{eqn:connection},
	\begin{align*}
	& \underset{A}{\inf}\|\Psi(\tau)-\Psi(\tau_0) \|_{\mathbb{L}\times \mathbb{R}^q}\\
	= & \underset{\{\gamma=0\} \cap \{\alpha,\gamma:\|\alpha-\alpha_0\|+\|\gamma-0\|\geq \epsilon\}}{\inf}\left(\| Q\psi_{1,\alpha,\gamma}^{c}-Q\psi_{1,\alpha_0,0}^{c} \|_H +\|Q \psi_{2,\gamma}^{c}-Q \psi_{2,0}^{c} \|_\mathbb{E}\right)\\
	=&\underset{\alpha:\|\alpha-\alpha_0\|\geq \epsilon}{\inf}\|Q\psi^*_{1,\alpha,0} -Q\psi^*_{1,\alpha_0,0}\|_H + 0\\
	=&\underset{\alpha:\|\alpha-\alpha_0\|\geq \epsilon}{\inf}\|P\psi_{\alpha} -P\psi_{\alpha_0}\|_H.
	\end{align*} 
	Under condition \circled{2}, the above display is greater than $0$.
	\begin{align*}
	&\underset{B}{\inf}\|\Psi(\tau)-\Psi(\tau_0) \|_{\mathbb{L}\times \mathbb{R}^q}\\
	= &\underset{B}{\inf}\left(\|Q\psi_{1,\alpha,\gamma}^{c}-Q\psi_{1,\alpha_0,0}^{c} \|_H + \|Q \psi_{2,\gamma}^{c}-Q \psi_{2,0}^{c} \|_\mathbb{E}\right)\\
	\geq & \underset{\{\gamma\neq 0\} \cap \{\alpha,\gamma:\|\alpha-\alpha_0\|+\|\gamma-0\|\geq \epsilon\}}{\inf} \|Q\psi^{c}_{2,\gamma}-Q \psi^{c}_{2,0}\|_\mathbb{E}\\
	\geq & \underset{\gamma:\|\gamma-0\| \geq \epsilon}{\inf} \|Q\psi^{c}_{2,\gamma}-Q \psi^{c}_{2,0}\|_\mathbb{E}.
	\end{align*}
	By lemma \ref{lemma:unique_solution}, $\gamma = 0$ is the unique solution to $Q\psi^{c}_{2,\gamma} = 0$. Because $\Gamma$ is compact and $Q\psi^{c}_{2,\gamma}(V,R)$ is a continuous function in $\gamma$, the unique solution $0$ implies that $\underset{\gamma: \|\gamma-0\| \geq \epsilon}{\inf}\|Q\psi^{c}_{2,\gamma}-Q \psi^{c}_{2,0}\|_\mathbb{E} > 0$ \citep{vdv1998}. %{\color{red}{Is there a better  reference here?}} 
	Therefore, 
	$\underset{B}{\inf}\|\Psi(\tau)-\Psi(\tau_0) \|_{\mathbb{L}\times \mathbb{R}^q} >0 $ and consequently
	$ \underset{A\cup B}{\inf}\|\Psi(\tau)-\Psi(\tau_0) \|_{\mathbb{L}\times \mathbb{R}^q}>0 $.
	%\note[Kate]{$ inf_{C} (A+B) =inf_{C}(A) +inf_{C}(B)$ because 1) for any $c \in C$, we obtain a pair of $(a,b)$. We always have $a \geq inf_C(A)$ and $b \geq inf_C(B)$, so $a+b \geq inf_{C} (A) +inf_{C}(B)$ and $inf_{C} (A+B) \geq inf_{C} (A) +inf_{C}(B)$ ; 2) Over the set C, there exists (a,b) such that $a < inf_{C} (A)+0.5\epsilon, b <inf_{C} (B)+0.5\epsilon$, thus $inf_{C}(A+B) \leq inf_{C} (A) + inf_{C} (B)$ } 
	%Because $\inf (D \cup E) = min \{\inf (D), \inf(E) \}$, we conclude for any $\epsilon>0$,
	% \begin{align*}
	% & \underset {(\alpha,\gamma):\|(\alpha,\gamma)-(\alpha_0,0)\|\geq \epsilon}{\inf}\|\Psi^{c}(\alpha,\gamma)-\Psi^{c}(\alpha_0,0) \|_{\mathbb{L}\times \mathbb{R}^q}\\
	% = & \underset{ A\cup B}{\inf} \|\Psi^{c}(\alpha,\gamma)-\Psi^{c}(\alpha_0,0) \|_{\mathbb{L}\times \mathbb{R}^q} \\
	% = & \min\left\{\underset{ A}{\inf} \|\Psi^{c}(\alpha,\gamma)-\Psi^{c}(\alpha_0,0) \|_{\mathbb{L}\times \mathbb{R}^q}, \underset{ B}{\inf} \|\Psi^{c}(\alpha,\gamma)-\Psi^{c}(\alpha_0,0) \|_{\mathbb{L}\times \mathbb{R}^q}\right \} >0.
	% \end{align*}
	% Condition \ref{condition2} is satisfied when we adapt it to the new map $\Psi^{c}$ and the new parameter $(\alpha, \gamma)$.
	Finally, by the definition of $\hat{\alpha}^{c}$ and $\hat{\gamma}$ in \eqref{EE:calibrated_estimators},  $\hat{\tau} = (\hat{\alpha}, \hat{\gamma})$ satisfies  the equation $\Psi_N(\hat{\tau})=0$.  It follows from theorem \ref{thm:consistency} that $\hat{\tau} \rightarrow_p \tau_0$ and we have $\hat{\alpha}^{c}\rightarrow_p \alpha_0$.
\end{proof}

Corollaries 1-3 have guaranteed our two-phase VPS estimator $\hat{\alpha}^*$ and the calibrated version of it $\hat{\alpha}^c$ will estimate the same target parameter $\alpha$ as the RS estimator $\hat{\alpha}$  as if we had obtained the complete information on $X$ for all phase I observations, regardless of the  model choice  for $X$ and the auxiliary variables choice $\tilde{V}$ for weights calibration.

Before we  provide the theoretical tools for deriving the limiting distributions of $\hat{\alpha}$, $\hat{\alpha}^*$ and $\hat{\alpha}^c$. We  introduce  the empirical process notation $\mathbb{G}_N$:
$
\mathbb{G}_N f= \sqrt{N}(\mathbb{P}_N f-Pf).
$
Suppose $f$ belongs to a class of functions $\mathcal{F}$,   $\mathbb{P}_N f$ then can be considered as a random map in $l^{\infty}(\mathcal{F})$ and the empirical process $\{\mathbb{G}_N f, f\in \mathcal{F}\}$ is the centered and scaled version of this random map. Each empirical process is associated with a class of functions $\mathcal{F}$. Therefore, each particular empirical process can be specified by an index set $\mathcal{F}$ and refereed as $\mathcal{F}$-indexed empirical process (see details in \cite{van2000weak}). In this paper, the $\mathcal{F}$-indexed empirical process we are concerned with can be both finite or infinite-dimensional depending on $\mathcal{F}$.

Let now discuss the asymptotic distribution of  $\hat{\alpha}$. With RS data,  Theorem 19.26 from \citet{vdv1998} can be immediately applied and we rephrase it as  Theorem \ref{thm:gauss:rs}.

\begin{theorem}(Aysmptotic normality of RS estimators)\label{thm:gauss:rs}
	%If $\|\mathbb{P}_N\psi_{\hat{\alpha}}\|_{ H}= o_p(N^{-1/2})$ and $\hat{\alpha} \overset{P}{\rightarrow} \alpha_0$,
	If $\psi_\alpha$ is a map satisfying conditions \circled{1} - \circled{5}, then
	\begin{equation}\label{eqn:normality_rs}
	\dot{\Psi}_0\sqrt{N}(\hat{\alpha} -\alpha_0) = - \mathbb{G}_N \psi_{\alpha_0} + o_p(1).
	\end{equation}
\end{theorem} 

\begin{proof}
	Under conditions \circled{1} - \circled{2}, it follows from corollary \ref{thm: consist-srs} that $\hat{\alpha} \rightarrow_p \alpha_0$.  According to the definition of $\hat{\alpha}$ in  \eqref{EE:complete}, $\mathbb{P}_N\psi_{\hat{\alpha}, h}=0, \forall h \in H$. These two results in addition to conditions \circled{3} - \circled{5} yield the desired results by theorem 19.26 \citep{vdv1998}.
\end{proof}
Let $\mathcal{F}_0 \equiv \{\psi_{\alpha_0, h}, h\in H \}$.  Then under condition \circled{3}, $\mathcal{F}_0$ is P-Donsker. By the definition of Donsker class \citep[Chap. 2.1]{van2000weak}, $\mathcal{F}_0$-indexed empirical process $\mathbb{G}_N\psi_{\alpha_0}$ in \eqref{eqn:normality_rs} converges in distribution to the $P$-Brownian bridge $\mathbb{G}$ in $l^{\infty}(\mathcal{F}_0)$. The limiting process $\{\mathbb{G}\psi_{\alpha_0,h}, \psi_{\alpha_0,h} \in \mathcal{F}_0\}$ is a zero-mean Gaussian process with covariance function
\[
E\mathbb{G}\psi_{\alpha_0,h}\mathbb{G}\psi_{\alpha_0,g}=P\psi_{\alpha_0,h}\psi_{\alpha_0,g} -P\psi_{\alpha_0,h}P\psi_{\alpha_0,g}
\]
where $h, g \in H$. Since $P\psi_{\alpha_0,h}=0$ and $P\psi_{\alpha_0,g}=0$, the asymptotic variance and covariance are reduced to
\begin{equation}\label{var:srs}
Var_A\left[\dot{\Psi}_0\sqrt{N}(\hat{\alpha} -\alpha_0)h\right] = P\psi_{\alpha_0,h}^2
\end{equation}
and 
\begin{equation}\label{cov:srs}
Cov_A(\dot{\Psi}_0\sqrt{N}(\hat{\alpha} -\alpha_0)h, \dot{\Psi}_0\sqrt{N}(\hat{\alpha} -\alpha_0)g)= P\psi_{\alpha_0,h}\psi_{\alpha_0,g}. 
\end{equation}

\begin{theorem}(Aysmptotic normality of two-phase VPS estimators)\label{thm:gauss:2ph} 
	Suppose assumptions A.\ref{assum: MAR} - \ref{assum: bounded_below} hold. If $\psi_\alpha$ is a map satisfying conditions \circled{1} - \circled{5} and it is uniformly bounded, then 
	\begin{equation}\label{eqn:normality_2ph}
	\dot{\Psi}_0\sqrt{N}(\hat{\alpha}^{*} -\alpha_0) = - \mathbb{G}_N \psi^*_{\alpha_0} + o_p(1).
	\end{equation}
\end{theorem}	
\begin{proof}
	We apply theorem 19.26 \citep{vdv1998} to $\psi^*_{\alpha}$ and $\hat{\alpha}^*$.  
		%     we have $\mathcal{G}' = \{\psi^*_{\alpha,h}(R,V,X) = \frac{R}{\pi_0(V)}\psi_{\alpha, h}(X): \alpha \in \mathbb{L}_0, h\in H\}$ is GC and condition \circled{2} holds when we replace $P\psi_{\alpha}, P\psi_{\alpha_0}$ by  $Q\psi_{\alpha}$ and $Q\psi_{\alpha_0}$.
	By assumption A  \ref{assum: bounded_below}, $\frac{R}{\pi_0(V)}$  is a uniformly bounded function. Under condition \circled{3},  $\mathcal{F}$ is Donsker and $\|P\psi_{\alpha} \|_{\mathcal{F}} < \infty$ for some $\delta  >0$.  It follows from \citet[example 2.10.10]{van2000weak} that $\mathcal{F}' \equiv \{\psi^*_{\alpha, h} =  \frac{R}{\pi_0(V)}\psi_{\alpha,h}: \|\alpha-\alpha_0\|< \delta, h \in H\}$ is Donsker for some $\delta>0$. This new class $\mathcal{F}' $ still has a finite and integrable envelope function as $\mathcal{F}$ does. 
	
	%According to this theorem, if we establish that $\{\psi^*_{\alpha,h},\|\alpha-\alpha_0\|<\delta, h\in H\}$ is Donsker for some $\delta >0$, with finite envelope function, that the map $\alpha \rightarrow Q\psi^*_{\alpha}$ is Fr\'echet differentiable at $\alpha_0$, with a derivative $\dot{\Psi}^*_0$ that has a continuous inverse on its range, and that $\|Q(\psi_{\alpha,h}^*-\psi_{\alpha_0,h}^*)^2 \|_{H} \rightarrow 0$ as $\alpha \rightarrow \alpha_0$, we assume $\hat{\alpha}^* \overset{P}{\rightarrow} \alpha_0$. then 
	% \[
	% \dot{\Psi}^*_0\sqrt{N}(\hat{\alpha}^{*} -\alpha_0) = - \mathbb{G}_N \psi^*_{\alpha_0} + o_p(1).
	%\]
	
	%First function $\psi^*_{\alpha,h}:(x,r,v) \rightarrow \psi^*_{\alpha, h} (r,x,v)$ is still measurable because $x \rightarrow \psi_{\alpha,h}(x)$ and $\frac{R}{\pi_0(V)}$ are both measurable, and measurable functions are closed under multiplication.
	
	%Under assumption \eqref{assum: bounded_below}, $\frac{R}{\pi_0(V)}$ is a fixed bounded function of random variables $R$ and $V$. Here ``fixed" means the function does not change over parameter spaces. It follows from Example 2.10.8 \citep{van2000weak} that $\{\psi_{\alpha,h}^*,\|\alpha-\alpha_0\|<\delta, h\in H\}$ as a product of a fixed bounded function and a Donsker class is also Donsker, with finite envelope function. 
In view of equation \eqref{eqn:connection}, the deterministic map $\alpha \rightarrow Q\psi^*_{\alpha}$ is the same as the map $\alpha \rightarrow P\psi_{\alpha}$. Thus, under condition \circled{4} the new map $\alpha \rightarrow Q\psi^*_{\alpha}$ is still Fr\'echet differentiable at $\alpha_0$, with a derivative $\dot{\Psi}_0$ that has a continuous inverse on its range. 
By A.\ref{assum: MAR},
	\begin{align*}
	Q(\psi_{\alpha,h}^*-\psi_{\alpha_0,h}^*)^2 
	=& Q\left[\frac{R^2}{\pi_0^2(V)}(\psi_{\alpha,h}-\psi_{\alpha_0,h})^2\right]\\
	=& 
	E_{Q}\left [\{\psi_{\alpha,h}(X)-\psi_{\alpha_0,h}(X)\}^2
	E_Q\left [ \frac{R}{\pi_0^2(V)}|X,U\right] \right] \\
	= & Q\left[(\psi_{\alpha,h}-\psi_{\alpha_0,h})^2\frac{1}{\pi_0^2(V)} 
	P(R=1|X,U) \right] \\
	= & Q\left[(\psi_{\alpha,h}-\psi_{\alpha_0,h})^2\frac{1}{\pi_0(V)}\right].
	\end{align*}
	Given  A.\ref{assum: bounded_below}, we have $\|Q(\psi_{\alpha}^*-\psi_{\alpha_0}^*)^2 \|_{H}  \leq  \frac{1}{\sigma}\|P(\psi_{\alpha}-\psi_{\alpha_0})^2 \|_{H} $. Under condition \circled{5}, $\|P(\psi_{\alpha}-\psi_{\alpha_0})^2 \|_{H} \rightarrow 0$ as $\alpha \rightarrow \alpha_0$.  Hence, $\|Q(\psi_{\alpha}^*-\psi_{\alpha_0}^*)^2 \|_{H} \rightarrow 0$  as $\alpha \rightarrow \alpha_0$.
	
	Based on  corollary \ref{thm: consist-2ph}, we  have $\hat{\alpha}^* \rightarrow_p \alpha_0$ and by \eqref{EE:twophase},  $\|\mathbb{P}_N\psi^*_{\hat{\alpha}^*}\|=0$  
	Applying theorem 19.26 \citep{vdv1998} to $\psi^*_{\alpha}$ we obtain $
	\dot{\Psi}_0\sqrt{N}(\hat{\alpha}^{*} -\alpha_0) = - \mathbb{G}_N \psi^*_{\alpha_0} + o_p(1).$
\end{proof}

Let $\mathcal{F}'_0 \equiv \{\psi^*_{\alpha_0, h}, h\in H\}$. Since we have shown  $\mathcal{F}'$ is a Donsker class with a finite and integrable envelope function for a small neighborhood of $\alpha_0$.  Thus, the $\mathcal{F}'_0$-indexed empirical process $\mathbb{G}_N\psi^*_{\alpha_0}$ in \eqref{eqn:normality_2ph} converges in distribution to a $Q$-Brownian bridge $\mathbb{G}$ in $l^{\infty}(\mathcal{F}'_0)$. 
%Given \eqref{eqn:connection} and \eqref{def:para}, $Q\left[\frac{R}{\pi_0(V)}\psi_{\alpha_0,h}\right] =0$ and $Q\left[\frac{R}{\pi_0(V)}\psi_{\alpha_0,g}\right]=0$. 
Applying the variance decomposition formula  we obtain 
\begin{align*}
Var\left[\psi^*_{\alpha_0,h}\right]=
& Var\left[E\left[\frac{R}{\pi_0(V)}\psi_{\alpha_0,h}|V\right] \right]
+E\left[Var\left[\frac{R}{\pi_0(V)}\psi_{\alpha_0,h}|V\right]\right]\nonumber\\
=& P\psi_{\alpha_0,h}^2+Q\left[\frac{1-\pi_0(V)}{\pi_0(V)}\psi^2_{\alpha_0,h}\right].
\end{align*}
%Following the the same derivation, we can also obtain the covariance functions in the space $l^{\infty}(\mathcal{F}')$ by substituting $\psi_{\alpha_0,h}\psi_{\alpha_0,g}$ for $\psi^2_{\alpha_0,h}.$

% is %simplified to 
%Comparing asymptotic variance for the improved two-phase VPS estimator in \eqref{var:2phase} to 
% Var_A\left[\dot{\Psi}_0\sqrt{N}(\hat{\alpha}^* -\alpha_0)h\right]\nonumber
%\[
%P\psi_{\alpha_0,h}\psi_{\alpha_0,g}+Q\left[\frac{1-\pi_0(V)}{\pi_0(V)}\psi_{\alpha_0,h}\psi_{\alpha_0,g}\right].
%\]
Therefore, the asymptotic variance 
\begin{equation}\label{var:2ph}
Var_A\left[\dot{\Psi}_0\sqrt{N}(\hat{\alpha}^* -\alpha_0)h\right]= P\psi_{\alpha_0,h}^2+Q\left[\frac{1-\pi_0(V)}{\pi_0(V)}\psi^2_{\alpha_0,h}\right].
\end{equation}
The asymptotic covariance has the same form as above except  that $\psi^2_{\alpha_0,h}$ is replaced by $\psi_{\alpha_0,h}\psi_{\alpha_0,g}$.
\begin{remark}
	The asymptotic variance in \eqref{var:2ph} has two components. The first equals  the variance we would obtain if 
	we had collected complete information on $X$ for all $N$ number of phase I subjects as shown in \eqref{var:srs}. %, in other words, if we had a random sample of $N$ subjects. 
	The second component is a penalty term for the fact that we only observe complete information on $X$ among the phase II subsample. 
\end{remark}

%is the additional variation introduced from phase II sampling, a penalty due to the fact that we do not observe complete information on $X$. 
% in reality. The second term also captures the difference between the asymptotic variance of $\hat{\alpha}$ and $\hat{\alpha}^*$ if we compare \eqref{var:2phase} to \eqref{var:srs}.

\begin{theorem}(Asymptotic normality of calibrated two-phase VPS estimators) \label{thm:gauss:cali}
	Suppose $ \Gamma$ is a compact convex subset of $\mathbb{R}^q $ with 0 as an interior point. Suppose  assumptions  A. \ref{assum: MAR}  - \ref{assum: vtilde-posdef} hold. If $\psi_\alpha$ is a map satisfying conditions \circled{1} - \circled{5} and it is uniformly bounded, then
	\begin{equation}\label{eqn:normality_cali}
	\dot{\Psi}^{c}_0\sqrt{N}\left(\begin{array}{c}
	\hat{\alpha}^{c} -\alpha_0\\
	\hat{\gamma}-0
	\end{array}
	\right) = - \mathbb{G}_N \left(\begin{array}{c}
	\psi^{c}_{1, \alpha_0,0} \\
	\psi^{c}_{2,0}
	\end{array} \right)+ o_p(1),
	\end{equation}
	where $\dot{\Psi}^{c}_0$ is the Fr\'echet derivative of the map $(Q \psi^{c}_{1,\alpha,\gamma}, Q \psi^{c}_{2,\gamma}): \mathbb{L}_0 \times \Gamma \mapsto \mathbb{L} \times \mathbb{R}^q$, defined by % $(\alpha, \gamma) \mapsto [\Psi^{c}_1(\alpha,\gamma),\Psi^{c}_2(\alpha,\gamma) ]$ at $(\alpha_0, 0)$. 
	% that has an continuous inverse on its range. $\dot{\Psi}^{c}_0 : lin \mathbb{L}_0 \times lin \Gamma \mapsto \mathbb{L} \times \mathbb{R}^q$ takes the form 
	% \[
	% (\alpha-\alpha_0, \gamma-\gamma_0)\mapsto \left(\begin{array}{cc}
	% \dot{\Psi}^{c}_{11} & \dot{\Psi}^{c}_{12}\\
	% \dot{\Psi}^{c}_{21} & \dot{\Psi}^{c}_{22}
	% \end{array}\right)\left(\begin{array}{c}
	% \alpha-\alpha_0\\
	% \gamma-\gamma_0
	% \end{array}\right)
	% \]
	\begin{equation}\label{eqn: frechet: cali}
	\dot{\Psi}^{c}_0 =\left(\begin{array}{cc}
	\dot{\Psi}^{c}_{11} & \dot{\Psi}^{c}_{12}\\
	\dot{\Psi}^{c}_{21} & \dot{\Psi}^{c}_{22}
	\end{array}\right)
	= \left(\begin{array}{cc}
	\dot{\Psi}_0 & 
	-Q\psi_{\alpha_0}\tilde{V}^T\\
	0 &
	-Q\tilde{V}\tilde{V}^T
	\end{array}\right).
	\end{equation}
	%Let $\mathcal{F}''= \{\psi^{c}_{1;\alpha,\gamma, h}, \alpha,\gamma \in \mathbb{L}_0 \times \Gamma, h \in H \}\cup \{\psi^{c}_{2,\gamma},\gamma \in \Gamma\}$, 
\end{theorem}

\begin{proof}
	Consider parameter $\tau = (\alpha, \gamma)$ and parameter space $\mathbb{L}_0 \times \Gamma$. Let $\psi_{\tau}(x, v, r)$ denote the map $c \mapsto \psi_{\tau, c}$ where $c = (h, b)$ and $\psi_{\tau, c}(x, v, r)= \psi^c_{1,\alpha,\gamma, h}(x, v, r) + b^T\psi^c_{2,\gamma}(v)$. We require $ c \in H \times B $  with $H$  an arbitrary index set and $B$ a unit ball in $\mathbb{R}^q$. We prove this theorem by applying theorem 19.26 \citep{vdv1998} to $\psi_{\tau}(x, v, r)$.

	%	Similar to the proof for theorem \ref{thm:gauss:2ph}, we adapt theorem \ref{thm:gauss:rs} to $(\hat{\alpha}^c, \hat{\gamma})$ and $(\psi^c_{1,\alpha,\gamma,h}, \psi^c_{2,\gamma})$ this time. Note $(\alpha, \gamma)$ can be identified as a map in $l^{\infty}(H \times B)$ that maps $(h, b)$  into $\mathbb{R}$ by $\alpha h + b^T\gamma$, where $B$ is a bounded subset of $\mathbb{R}^q$. Hence theorem \ref{thm:gauss:rs} is still  applicable.

	First we show  $\{\psi_{\tau, c} :  \|\tau-\tau_0\|_{\mathbb{L}_0 \times \Gamma} < \delta, c \in H \times B \}$ is Donsker with finite envelope function for some $\delta >0$.    $\{f_{\gamma},\gamma \in \Gamma\}$ is uniformly bounded and uniformly Lipschitz. It follows from example 19.7 \citep{vdv1998} that $\{f_{\gamma}, \gamma \in \Gamma\}$ is Donsker.  Under A.\ref{assum: vtilde-bounded}, a singleton $\{\tilde{V}\}$ is a finite set of square integrable functions, thus it is Donsker. Since $b^T\tilde{V}$ is a Lipschitz function, $\{b^T\tilde{V}, b \in B \}$  is Donsker as well.  In addition, condition \circled{3} ensures $\{\psi_{\alpha,h}: \|\alpha-\alpha_0\|< \delta, h \in H\}$ is  Donsker for some $\delta >0$. By assumption, functions in these three Donsker classes are all uniformly bounded. Given $\psi_{\tau, c} = f_{\gamma} \psi_{\alpha, h} +  f_{\gamma}b^T\tilde{V}-b^T\tilde{V}$,  it follows from the preservation theorem of Donsker\citep{van2000weak} that $\{\psi_{\tau, c} :  \|\tau-\tau_0\|_{\mathbb{L}_0 \times \Gamma} < \delta, c \in H \times B \}$ is Donsker with a finite envelope function for some $\delta >0$.

$\dot{\Psi}^{c}_0$ given in \eqref{eqn: frechet: cali}  can be easily obtained by following the definition of Fr\'echet derivative. Details are provided in Appendix B. To ascertain $\dot{\Psi}^{c}_0$ has a continuous inverse  on its range, we only need to verify that $\dot{\Psi}^{c}_{11}$ and $\dot{\Psi}^{c}_{22}-\dot{\Psi}^{c}_{21}\dot{\Psi}^{**^{-1}}_{11}\dot{\Psi}^{c}_{12}$ have continuous inverses  \citep[ch. 25]{vdv1998}. Because $\dot{\Psi}^{c}_{11}=\dot{\Psi}_0$, condition \circled{4} ensures  that $\dot{\Psi}^{c}_{11}$ has a continuous inverse. Since $\dot{\Psi}^{c}_{21}=0$, $\dot{\Psi}^{c}_{22}-\dot{\Psi}^{c}_{21}\dot{\Psi}^{**^{-1}}_{11}\dot{\Psi}^{c}_{12}= \dot{\Psi}^{c}_{22} = Q\tilde{V}\tilde{V}^T$. Under assumption A.\ref{assum: vtilde-posdef},  $Q\tilde{V}\tilde{V}^T$ is positive definite and thus  $Q\tilde{V}\tilde{V}^T$ has an inverse. Since $Q\tilde{V}\tilde{V}^T$ is a linear maps between Euclidean spaces,  it is also continuous.  
	
%	{\color{red}{The old proof $Q\tilde{V}\tilde{V}^T$ is continuous. Under assumption A.\ref{assum: vtilde-posdef}, the null space of $Q\tilde{V}\tilde{V}^T$ is $0$. By assumption A.\ref{assum: vtilde-bounded}, $\tilde{V}$ is also bounded. Since the domain of $Q\tilde{V}\tilde{V}^T$ is Lin$\Gamma$ and it is closed, thus the range of $Q\tilde{V}\tilde{V}^T$ is closed. According to \citet [proposition 7B in Appendix]{bickel1993efficient}, $Q\tilde{V}\tilde{V}^T$ has an continuous inverse as well. }}
	%Second, we adapt condition \ref{con:frechet} to $\psi^{c}_{\tau,c}$ and show as a map into $l^\infty (\mathbb{C}),$ the map $\tau \mapsto \Psi^{c}_{\tau}$ is Fr\'echet-differentiable at $\tau_0 = (\gamma_0,0)$, with a derivative $\dot{\Psi} ^{c}_0: lin \mathbb{T} \mapsto l^\infty ( \mathbb{C})$ that has a continuous inverse on its range. We prove the existence of $\dot{\Psi} ^{c}_0$ by validating each component in its explicit form 
	%given in \eqref{eqn: frechet: cali}. 
	
	Next,we verify  $\|(\psi_{\tau, c}-\psi_{\tau_0, c})^2\|_{H \times B} \rightarrow 0 $ as 	$\tau \rightarrow \tau_0$.
	\begin{align*}
	&\sup_{h \in H} Q(\psi_{1,\alpha,\gamma, h}^{c}-\psi_{1,\alpha_0,\gamma_0, h}^{c})^2\\
	= & \sup_{h \in H} Q(\psi_{1,\alpha,\gamma, h}^{c}-\psi_{1,\alpha_0,\gamma, h}^{c} +\psi_{1,\alpha_0,\gamma, h}^{c} -\psi_{1,\alpha_0,\gamma_0,h}^{c})^2 \\
	\leq & \sup_{h \in H} 2Q(\psi_{1,\alpha,\gamma, h}^{c}-\psi_{1,\alpha_0,\gamma, h}^{c})^2 +\sup_{h \in H} 2Q(\psi_{1,\alpha_0,\gamma, h}^{c} -\psi_{1,\alpha_0,\gamma_0,h}^{c})^2\\
	= & 2\sup_{h \in H}Q\left[\frac{R^2}{\pi_0(V)^2}\exp(-\gamma^T\tilde{V})^2(\psi_{\alpha,h}-\psi_{\alpha_0,h})^2\right]\\
	& \qquad + 2\sup_{h \in H}Q\left[\frac{R^2}{\pi_0(V)^2}\psi_{\alpha_0,h}^2\left\{\exp(-\gamma^T\tilde{V})-\exp(-\gamma_0^T\tilde{V})\right\}^2\right].
	\end{align*}
	
 Based on our assumptions, $R, 1/\pi_0(V),\gamma$,  and $\tilde{V}$ are all bounded and $\psi_{\alpha,h}(x)$ is uniformly bounded. Thus, there exist some positive constants $C_1$ and $C_2$ such that
	\begin{align*}
	&\sup_{h \in H} Q(\psi_{1,\alpha,\gamma, h}^{c}-\psi_{1,\alpha_0,\gamma_0, h}^{c})^2\\
	\leq & C_1\sup_{h \in H}|P(\psi_{\alpha,h}-\psi_{\alpha_0,h})^2|+C_2Q\left[\exp(-\gamma^T\tilde{V})-\exp(-\gamma_0^T\tilde{V})\right]^2.
	\end{align*}
	Under condition \circled{5}, the first term on the RHS  of the above inequality goes to 0 as $\alpha \rightarrow \alpha_0$. By dominated convergence theorem, the second term also goes to 0 as $\gamma \rightarrow \gamma_0$.  Thus,   $
	\sup_{h \in H} Q(\psi_{1,\alpha,\gamma, h}^{c}-\psi_{1,\alpha_0,\gamma_0, h}^{c})^2 \rightarrow 0\mbox{ as } (\alpha, \gamma) \rightarrow (\alpha_0, \gamma_0)$. Using a similar argument, we can show  $b^TQ(\psi_{2,\gamma}^{c}-\psi_{2,\gamma_0}^{c})^2 \rightarrow 0$ as $\gamma \rightarrow \gamma_0$. Since $f_{\gamma}$ and $\psi_{\alpha}$ are both uniformly bounded, and $\psi^c_{1,\alpha,\gamma, h} = f_{\gamma}\psi_{\alpha}$, 
	$(\psi^c_{1,\alpha,\gamma, h}-\psi^c_{1,\alpha_0,\gamma_0, h})$ is uniformly bounded as well.  Using the dominance convergence theorem again we have $2 \sup_{h \in H, b \in B}Q|(
	\psi^c_{1,\alpha,\gamma, h}-\psi^c_{1,\alpha_0,\gamma_0, h})b^T(\psi_{2, \gamma}^{c}-\psi^c_{2,\gamma_0})| \rightarrow 0 $ as $(\alpha, \gamma) \rightarrow (\alpha_0, \gamma_0)$.
	
Because
	\begin{align*}
	&\|(\psi_{\tau, c}-\psi_{\tau_0, c})^2\|_{H \times B}\\
	= & \sup_{h \in H, b\in B} Q(
	\psi_{1,\alpha,\gamma, h}^{c} -\psi_{1,\alpha_0,\gamma_0,h}^{c} + 
	b^T\psi_{2,\gamma}^{c} -b^T\psi_{2, \gamma_0}^{c} )^2 \\
	\leq & \sup_{h \in H} Q(
	\psi_{1,\alpha,\gamma, h}^{c}-\psi_{1,\alpha_0,\gamma_0, h}^{c})^2 +
	b^TQ( \psi_{2, \gamma}^{c} -\psi_{2,\gamma_0}^{c})^2 \\
	&\qquad  +2 \sup_{h \in H, b\in B}Q|(
	\psi^c_{1,\alpha,\gamma, h}-\psi^c_{1,\alpha_0,\gamma_0, h})b^T(\psi_{2, \gamma}^{c}-\psi^c_{2,\gamma_0})|,
	\end{align*}
 and we have shown that the three terms on the RHS of the above display all converge to $0$ as $(\alpha, \gamma) \rightarrow (\alpha_0, \gamma_0)$, condition \circled{5} is satisfied for $\psi_{\tau}$ in $l^{\infty}(H \times B)$. 
	
	Finally, by corollary \ref{thm: consist-cali}, $\hat{\tau} \rightarrow_{p} \tau_0$  and  by the definition  of $(\hat{\alpha}^{c},\hat{\gamma})$ in \eqref{EE:calibrated_estimators}, $\mathbb{P}_N\psi_{\hat{\tau}}(X, V, R) = 0$.  Applying theorem 19.26 \citep{wellner1992} to $\psi_{\tau}$ we obtain 
	$\dot{\Psi}^{c}_0\sqrt{N}\left(\begin{array}{c}
	\hat{\tau}-\tau_0
	\end{array}
	\right) = - \mathbb{G}_N \psi_{\tau_0} + o_p(1)$.
\end{proof}
The display in  \eqref{eqn:normality_cali} is equivalent to for all $h \in H$, 
\begin{align*}
\sqrt{N}\dot{\Psi}_0 (\hat{\alpha}^c-\alpha_0)h -\sqrt{N}Q [\psi_{\alpha_0,h}(X)\tilde{V}^T](\hat{\gamma}-\gamma_0) & = - \mathbb{G}_N\psi^{c}_{\alpha_0,h}+ o_p(1)\\
-\sqrt{N}Q\tilde{V}\tilde{V}^T(\hat{\gamma}-\gamma_0) & = - \mathbb{G}_N\psi^{c}_{2,0} + o_p(1).
\end{align*}
To find the limiting distribution of $\hat{\alpha}^{c}$ alone, we cancel out the terms with $(\hat{\gamma} -\gamma_{0})$ by applying the linear operator $Q[\psi_{\alpha_0,h}(X)\tilde{V}^T](Q\tilde{V}\tilde{V}^T)^{-1}$ to both sides of the first equation and then subtract it from the second equation. As a result,
\begin{align*}
&\sqrt{N}\dot{\Psi}_0 (\hat{\alpha}^{c}-\alpha_0)h \\
= & -\mathbb{G}_N\psi^{*}_{\alpha_0,h}+Q\left[\psi_{\alpha_0,h}(X)\tilde{V}^T\right](Q\tilde{V}\tilde{V}^T)^{-1}\mathbb{G}_N\psi^{c}_{2,0}+o_p(1)\nonumber\\
= & -\mathbb{G}_N\frac{R}{\pi_0(V)}\psi_{\alpha_0,h}+\mathbb{G}_N\{Q\psi_{\alpha_0,h}(X)\tilde{V}^T\}(Q\tilde{V}\tilde{V}^T)^{-1}\{\frac{R}{\pi_0(V)}-1\}\tilde{V}+o_p(1)\nonumber\\
%= & -\mathbb{G}_N\psi_{\alpha_0,h}-\mathbb{G}_N\{\frac{R}{\pi_0(V)}-1\} \left[\psi_{\alpha_0,h}-Q\{\psi_{\alpha_0,h}(X)\tilde{V}^T\}(Q\tilde{V}\tilde{V}^T)^{-1}\tilde{V}\right]+o_p(1)\\
= & -\mathbb{G}_N\left[\psi_{\alpha_0,h}+\{\frac{R}{\pi_0(V)}-1\}\left\{\psi_{\alpha_0,h}-\Pi(\psi_{\alpha_0,h}|\tilde{V})\right\}\right]+o_p(1)
\end{align*}
where $\Pi(\cdot|\tilde{V})$ denotes the population least squares projection on the space spanned by the calibration variable $\tilde{V}: \Pi(\cdot|\tilde{V})=Q\{\cdot\tilde{V}^T\}(Q\tilde{V}\tilde{V}^T)^{-1}\tilde{V}.$

Therefore, $\sqrt{N}\dot{\Psi}_0 (\hat{\alpha}^{c}-\alpha_0)$ %is asymptoticly equivalent to an empirical process indexed $\mathbb{G}$ by $\mathcal{F} \times \mathcal{F}_4$. This empirical process
converges in distribution to the Q-Brownian bridge in $l^{\infty}(\mathcal{F}_0'')$ where $\mathcal{F}_0'' = \left\{\psi_{\alpha_0,h}+\{\frac{R}{\pi_0(V)}-1\}\{\psi_{\alpha_0,h}-\Pi(\psi_{\alpha_0,h}|\tilde{V})\}, h \in H \right\} $
This Q-Brownian bridge is a zero-mean Gaussian with variance 
\begin{equation}\label{var:cali}
Var_A\{\sqrt{N}\dot{\Psi}_0 (\hat{\alpha}^{c}-\alpha_0)h\}
= P\psi^2_{\alpha_0,h} +Q\left[\frac{1-\pi_0(V)}{\pi_0(V)}\left\{\psi_{\alpha_0,h}-\Pi(\psi_{\alpha_0,h}|\tilde{V})\right\}^2\right].
\end{equation}
Replacing $\psi^2_{\alpha_0, h}$ in \eqref{var:cali} with $\psi_{\alpha_0, h}\psi_{\alpha_0, g}$ yields the  asymptotic covariance.
\begin{remark}
	Comparing the asymptotic variance for the calibrated two-phase VPS estimator in \eqref{var:cali} to the two-phase VPS estimator in \eqref{var:2ph}, we find that the auxiliary information will reduce the variance of a two-phase VPS  estimator if 
	\[
	Q\left[\frac{1-\pi_0(V)}{\pi_0(V)}\left\{\psi_{\alpha_0,h}-\Pi(\psi_{\alpha_0,h}|\tilde{V})\right\}^ 2\right]\leq Q\left[\frac{1-\pi_0(V)}{\pi_0(V)}\psi_{\alpha_0,h}^2\right].
	\]
	%where $ h \in H$.
	%In addition,
	%$\sqrt{N}\dot{\Psi}_0 (\hat{\alpha}^{c}-\alpha_0)$ converges to a Q-Brownian bridge in $l^{\infty}(H)$ with mean 0 and 
	%asymptotic variance 
\end{remark}

	The asymptotic variances results in \eqref{var:srs},  \eqref{var:2ph} and \eqref{var:cali} are all robust variances since we estimate allowing $P \in \mathcal{P}$. When $P \in\{ P_{\alpha}, \alpha \in \mathbb{L}_0\}$, by simply replacing $P$ in these results with $P_{\alpha}$ we obtain the model-based variances. Based on our proofs, our Z-estimators do not need to be the exact solutions to the EE. The theoretical results still hold for $\hat{\alpha}$ if  $\mathbb{P}_N\psi_{\hat{\alpha}} = o_p(1)$ and similarly for $\hat{\alpha}^*$ and $\hat{\alpha}^c$.

\section{An example}\label{sec:example}
 We illustrate the use of the Z-estimation system by applying it to Lin and Ying additive hazards (AH) model \citep{lin1994}. This model specifies the hazard function of a censored failure time $T$ as a sum of a baseline hazard function $\lambda(\cdot)$ and a regression function of $Z$:
\begin{equation}\label{model:LY}
\lambda (t|Z) = \lambda(t)+Z^T\theta.
\end{equation}

We will show how we can systematically develop estimators from the three data types we have discussed in section \ref{sec:prob} for three parameters --- the regression coefficients $\theta$, the cumulative baseline hazard $\Lambda(\cdot)=\int_0^{\cdot}\lambda(t)dt$, and the individual risk $\Lambda(t|Z = z)$ at time $t$ based on an individual's characteristics $z$. The first two parameters are of great interest in public health and the third one is useful for personalized risk prediction. 

In the following, we will provide key steps for developing these estimators. We refer interested readers to \citet{kate2014} for more detailed steps and numerical validation of our theoretical results. We first describe the data we consider and the assumptions we make. Let $\tilde{T}$ denote the true failure time, $C$ the censoring time, and $Z$ the p-dimensional covariates. We assume that there exists a finite maximum censoring time $\tau$ such that $Pr(C\geq \tau) = Pr(C=\tau) > 0$ and  $\tilde{T}$ is independent of $C$  given $Z$. Then $T = \tilde{T} \wedge C$ is the censored failure time and 
$\Delta = 1(\tilde{T} \leq C)$ is the indicator of whether we observe the true failure time.  The information of main interest in this problem can be summarized  as $X = (T, \Delta, Z)$. For RS data, $X$ are observed for all the samples.  With two-phase VPS, some covariates in $Z$ are only measured for the selected phase II subsample and thus only these members have complete information of $Z$ and therefore $X$.  In addition, we will also observe some auxiliary variables $U$  for all the phase I sample.

%When a two-phase VPS is implemented, a part of $X$ are, however, missing for the phase I members who are not selected into the phase II subsample. In the meantime, it is likely we also observe some auxiliary variables $U$  for all the phase I members. For example, in a  biomarker epidemiology study, the entire cohort (the phase I members) are followed for the outcomes $(T, C)$. In addition, some of baseline covariates of interest and some auxiliary variables $U$ are measured. Based on these covariates, the phase II subsample is determined and furthered measured for an expensive biomarker. Only the phase II subsample members have complete information of the covariates of interest $Z$ and therefore $X$.

 Our notations now follow  those used in section \ref{sec:sampling}. We make the same assumptions as A.\ref{assum: i.i.d.} - \ref{assum: vtilde-posdef} for our data and sampling scheme. We also assume 
 \begin{assumption}\label{assum: boundedness} 
 $Z$ is bounded and $PY(\tau)$ is bounded from 0.
 \end{assumption}
 \begin{assumption} \label{assum: A-posdef}
  \begin{equation}\label{eqn:A}
A \equiv P\int_0^\tau[Z-\frac{PZY(t)}{PY(t)}]^{\otimes 2}Y(t)dt
\end{equation}
is positive definite. 
\end{assumption}
Next, we follow the estimating procedures introduced in section \ref{sec:Z-system} to develop our estimators in four steps.
\subsection*{Step 1: Define Parameters} Instead of estimating the three parameters $\theta$, $\Lambda$, and $\Lambda(t|Z=z)$  individually, we  treat $\theta$ and $\Lambda$  as a single parameter and estimate them jointly. 
Let $\mathbb{L}_0$ be the parameter space for $\alpha =(\theta, \Lambda)$. Consider the following two  index sets
\begin{align}\label{function: h}
H_1 & \equiv \{ h_1 \in \mathbb{R}^p, \|h_1\|_E \leq M\} \mbox{ and } \\
H_2  & \equiv \{ h_2 \in BV[0,\tau], \|h_2\|_{H_2}\leq M, \|h_2\|_{TV} \leq M\}
\end{align}
where $\|\cdot\|_{TV}$ is a total variation norm and $ M $ is a finite positive constant. %i.e., $\|h_2\|_{TV} = \sup\sum_{i=1}^n |h_2(x_{i})-h_2(x_{i-1})|$ where the supreme is taken over all partitions $0 = x_0 < x_1, \dots, x_{n-1} < x_n =\tau$. %$H_1$ is the unit ball of the Euclidean space $\mathbb{R}^p$. $H_2$ is the set of all uniformly bounded functions of uniformly bounded variation over $[0, \tau]$, with both bounds equal to 1. 
%Let $\mathbb{L}_0$ be the parameter space for $(\theta, \Lambda)$.  Conventionally, $\mathbb{L}_0$ is considered as a product space $\Theta \times \mathbb{A} $ where $\Theta$ is a bounded subset of $\mathbb{R}^p$ and $\mathbb{A}$ is a collection of finite nondecreasing and nonnegative functions over the time interval $[0,\tau]$.
 Let $H  \equiv H_1 \times H_2$. Then the joint parameter $(\theta, \Lambda)$ can be uniquely identified by the map $ \alpha: H \mapsto\mathbb{R}$ defined as $\alpha h = h_1^T\theta+\int_0^\tau h_2d\Lambda$. Therefore, $(\theta, \Lambda) \in l^{\infty}(H)$ and  $(\theta, \Lambda)$ belongs to the type of parameter we consider in our Z-estimation system. 
 % \begin{remark} 
 % Our restriction on the index set $H$  is mostly for technical reasons. Later we will see our choice on $H$ facilitates us to establish asymptotic normality of our estimators.
 %  Other potentially useful choices of $H$ are discussed in \citet{murphy2001semiparametric}. 
 %  The nontechnical reason is that we make $H$ large enough to include $\Lambda$ in the usual sense of a cumulative hazard function over $[0,\tau]$. This is satisfied by our choice of $  H$ since $\forall s \in [0,\tau]$, $\Lambda(s)$ can be written as $\int 1_{[0,s]}(t)d\Lambda(t)$, if we let $h_2= 1_{[0,s]}(t)$,  then $h_2 \in H_2$ and $\Lambda h_2 = \Lambda(s)$.  
 % \end{remark}

In this system the parameter is defined through a zero-value function $\psi_{\alpha,h}$.  In this example, 
let 
$N(t)= 1(T\leq t,\Delta =1)$ be a counting process and $Y(t)= \textbf{1}[T\geq t]$ be an at-risk process.
We define
\begin{equation}\label{function: psi}
\psi_{\alpha, h} = \psi_{1,\theta, \Lambda,h_1} + \psi_{2,\theta, \Lambda,h_2}
\end{equation}
where
\begin{align*}
\psi_{1,\theta, \Lambda,h_1}(X)=& h_1^T\int_0^{\tau} Z\{dN(t) - Y(t) d\Lambda(t) - Y(t)Z^T\theta dt\},\\
\psi_{2,\theta, \Lambda,h_2}(X)= & \int_0^{\tau} h_2(t)\{dN(t)-  Y(t)d\Lambda(t) - Y(t)Z^T\theta dt\}.
\end{align*}
 Our $\psi_{\theta, \Lambda, h}$ is motivated from  the tangent set of the assumed AH model $P_{\theta, \Lambda}$ with respect to $\theta$ and $\Lambda$ in \eqref{model:LY}.  The detailed  derivation of our   $\psi_{\theta, \Lambda, h}$ is provided in \citet[ch. 3.6]{kate2014}.

Following section \ref{sec:para}, we define the true parameter $(\theta_0, \Lambda_0)$ as the solution to
\begin{equation}\label{eqn:para:ah}
 \Psi(\theta,\Lambda)h = P\psi_{\Lambda,\theta,h}(X) = 0,  \forall h \in H.
\end{equation}
Although  $\psi_{\theta, \Lambda,h}$ is motivated from $P_{\theta, \Lambda}$, we allow the true model $P$ to be any model in $\mathcal{P}$ with $PXX^T < \infty$.

% It can be shown $(\theta_0, \Lambda_0)$ is an explict solution to equation \eqref{eqn:para:ah}, so it is unique. Thus $(\theta_0, \Lambda_0)$ is well defined. We further add assumptions that $Z$ is bounded and $PY(\tau)$ is bounded from 0, so that $\psi_{\theta, \Lambda,h}$ also satisfies the requirement that $\psi_{\theta, \Lambda,h}$  and $P\psi_{\theta, \Lambda,h}$ are uniformly bounded.

\subsection*{Step 2: Construct Estimators} 
 Following  section \ref{sec:estimator},  we construct our EE based on  \eqref{EE:complete},  \eqref{EE:twophase},  or \eqref{EE:calibrated_estimators} for each data type respectively. Our estimators are obtained by solving these EE. The procedure to solve for our estimators can be  very systematic. Take RS data as an example. %to solve for $(\hat{\theta}, \hat{\Lambda})$,% we make use of the fact that EE \eqref{eqn:ee} holds for every $h \in H$ and therefore for every $h \in H\cdot M $ where $M$ is any finite positive number.   Thus we can first choose %The fact that our EE \eqref{eqn:ee}  holds for every $h\in H$ gives us lots of convenience in deriving our estimators. We first select
 We first chose an $h \in H$ that removes  $\Lambda$ and solve  for $\hat{\theta}$, and then we choose a different  $h \in H$ to obtain $\hat{\Lambda}$. For clarity, we suppress the notation $t$ in $ Y(t)$ in the middle of the derivation and put it back in the result. We set
\begin{equation}\label{fun: h-for-theta}
h= (h_1, h_2)= (h_1, -h_1^T\frac{\mathbb{P}_NZY}{\mathbb{P}_NY})
\end{equation}
 where $h_1$ is given in \eqref{function: h}. % For a particular dataset $\frac{\mathbb{P}_NZY}{\mathbb{P}_NY}$ is fixed, so $h_2 = h_1^T\frac{\mathbb{P}_NZY}{\mathbb{P}_NY}$ is valid. 
With this $h$, the LHS of EE \eqref{EE:complete} becomes 
\begin{align*}
\mathbb{P}_N\psi_{\theta, \Lambda,h} = &  \mathbb{P}_N(\psi_{1,\theta, \Lambda,h_1}+ \psi_{2,\theta, \Lambda,h_2})\\
= & \int_0^{\tau}\left\{ h_1^T\mathbb{P}_N ZdN(t) - h_1^T\mathbb{P}_N ZY d\Lambda(t) - h_1^T\mathbb{P}_N ZZ^TY\theta dt\right\}\\
& -\int_0^{\tau}\left\{
h_1^T\frac{\mathbb{P}_NZY}{\mathbb{P}_NY}\mathbb{P}_NdN(t)- h_1^T\frac{\mathbb{P}_NZY}{\mathbb{P}_NY}\mathbb{P}_NY d\Lambda(t) - h_1^T\frac{\mathbb{P}_NZY}{\mathbb{P}_NY} \mathbb{P}_NZ^TY\theta dt\right\} \\
= & h_1^T\int_0^{\tau}\mathbb{P}_N \left(Z-\frac{\mathbb{P}_NZY}{\mathbb{P}_NY}\right)dN(t)-h_1^T\int_0^\tau \mathbb{P}_N \left(Z-\frac{\mathbb{P}_NZY}{\mathbb{P}_NY}\right)Z^TY\theta dt.
\end{align*}
Because
\[
 \int_0^{\tau}\mathbb{P}_N \left(Z-\frac{\mathbb{P}_NZY}{\mathbb{P}_NY}\right)\frac{\mathbb{P}_NZ^TY}{\mathbb{P}_NY}Y\theta dt
=   \int_0^{\tau} \left[\mathbb{P}_N ZY\frac{\mathbb{P}_NZ^TY}{\mathbb{P}_NY}-\frac{\mathbb{P}_NZY}{\mathbb{P}_NY}\mathbb{P}_NZ^TY\right]\theta dt=0,
\]
Adding this zero-valued function to $\mathbb{P}_N\psi_{\theta, \Lambda,h}$ yields 
% $\int_0^{\tau}\mathbb{P}_N \left(Z-\frac{\mathbb{P}_NZY}{\mathbb{P}_NY}\right)\frac{\mathbb{P}_NZ^TY}{\mathbb{P}_NY}Y\theta dt $ to  both sides of EE, we obtain $\hat{\theta}$ solves
\[
\mathbb{P}_N\psi_{\theta, \Lambda,h} =  \mathbb{P}_N\int_0^{\tau} h_1^T\left(Z-\frac{\mathbb{P}_NZY}{\mathbb{P}_NY}\right)dN(t)-\mathbb{P}_N\int_0^\tau h_1^T\left(Z-\frac{\mathbb{P}_NZY}{\mathbb{P}_NY}\right)\left(Z-\frac{\mathbb{P}_NZY}{\mathbb{P}_NY}\right)^TY\theta dt.
\]
By definition, our estimator $\hat{\theta}$ solves $\mathbb{P}_N\psi_{\theta, \Lambda,h} = 0$. Rewriting the equation in a vector form, we have 
\[
\mathbb{P}_N\int_0^{\tau} \left(Z-\frac{\mathbb{P}_NZY}{\mathbb{P}_NY}\right)dN(t)-\left[Z-\frac{\mathbb{P}_NZY}{\mathbb{P}_NY}\right]^{\otimes 2}Ydt =0.
\]

$\mathbb{P}_N\int_0^{\tau}\left[Z-\frac{\mathbb{P}_NZY}{\mathbb{P}_NY}\right]^{\otimes 2}Ydt$ converges to A defined in \eqref{eqn:A} as $N \rightarrow \infty$. Under the assumption of A. \ref{assum: A-posdef}, when $N$ is large enough, $\mathbb{P}_N\int_0^{\tau}\left[Z-\frac{\mathbb{P}_NZY}{\mathbb{P}_NY}\right]^{\otimes 2}Ydt$ is positive definite and has an inverse. Therefore, 
\begin{equation}\label{theta_hat:RS}
\hat{\theta} = \left[\mathbb{P}_N\int_0^{\tau}\left\{Z-\frac{\mathbb{P}_NZY(t)}{\mathbb{P}_NY(t)}\right\}^{\otimes 2}Y(t)dt\right]^{-1}\mathbb{P}_N\int_0^{\tau}\left\{Z-\frac{\mathbb{P}_NZY(t)}{\mathbb{P}_NY(t)}\right\}dN(t).
\end{equation}\\
Next given we set 
\begin{equation}\label{fun: h-for-Lambda}
h=(h_1,h_2)=(0, \frac{1(t \leq s)}{\mathbb{P}_NY})
\end{equation}
where  $s \in [0, \tau]$.
Plugging $h$ into EE \eqref{EE:complete} yields
\[
\mathbb{P}_N\psi_{\theta, \Lambda,h} = 0 +
\int_0 ^{\tau}\left[ \frac{1(t\leq s)}{\mathbb{P}_NY}\mathbb{P}_NYdN(t)-\frac{1(t\leq s)}{\mathbb{P}_NY}\mathbb{P}_NYd\Lambda(t)-\frac{1(t\leq s)}{\mathbb{P}_NY}\mathbb{P}_NZ^TY\theta dt\right]=0.
\]
Replacing $\theta$ in the above equation by  $\hat{\theta}$ in \eqref{theta_hat:RS}, we obtain
\begin{equation}\label{Lambda_hat:RS}
\hat{\Lambda}(s)=\int_0^s\frac{\mathbb{P}_N[Y(t)dN(t)]}{\mathbb{P}_NY(t)} -\int_0^s\frac{\mathbb{P}_N[Z^TY(t)]}{\mathbb{P}_NY(t)}\hat{\theta}dt.
\end{equation}
Consequently, $\hat{\Lambda}(t|Z=z_)= \hat{\Lambda} + z^T\hat{\theta}$. Although this approach is  different from the one used by  \citet{lin1994} for developing their RS estimators, our results  are the same as theirs. Using a similar approach, we can derive $(\hat{\theta}^*, \hat{\Lambda}^*)$ and $(\hat{\theta}^c, \hat{\Lambda}^c)$ very fast, since a similar $h$ to \eqref{fun: h-for-theta} or \eqref{fun: h-for-Lambda} can be reused for the new data types. 

\subsection*{Step 3: Verify Conditions} Given $\Lambda$ interpreted as the  cumulative hazard function and $\theta$ as the regression coefficients, it is reasonable to assume
 \begin{assumption}\label{assum:compact}
 	The parameter space $\mathbb{L}_0$ for $(\theta, \Lambda)$ is compact, closed and bounded.  $\Lambda$ is a uniformly bounded function of uniformly bounded variation over $[0,\tau]$. 
 	\end{assumption}
 	%{\color{red}{In there a better way to state this assumption? It sounds repetitive to me.}}

Under our assumptions A \ref{assum: boundedness} - \ref{assum:compact}, 
we sketch our verification  of the five  conditions listed in  section \ref{sec:conditions}.
We first show condition \circled{2} that $(\theta_0, \Lambda_0)$ is a unique and well-separated solution to the equation  \eqref{eqn:para:ah}. Replacing $\mathbb{P}_N$ with $P$ in \eqref{fun: h-for-theta} and \eqref{fun: h-for-Lambda} generates two new $h$ indexes. Plugging them  into \eqref{eqn:para:ah} and following the same procedure  in Step 2 we will obtain  $(\theta_0, \Lambda_0)$.  $(\theta_0, \Lambda_0)$ is in the same form as \eqref{theta_hat:RS} and \eqref{Lambda_hat:RS} with $\mathbb{P}_N$ replaced by $P$.  $(\theta_0, \Lambda_0)$ as an explicit solution to \eqref{eqn:para:ah} is automatically unique. Because $\mathbb{L}_0$ is compact, the inequality in condition \circled{2} holds and $(\theta_0, \Lambda_0)$ is also a well-separated solution.% {\color{red}{$<$---any good reference here?}} 
Next, we verify conditions \circled{1} and \circled{3}. This can be proved by first  showing  each component of $\psi_{\theta, \Lambda, h}$ forms a Donsker class. These components include $Z$, $\Delta$,  $T$, $Z^T\theta$, $\Lambda(T)$, $h_2(T)$, $\int_0^\tau h_2(t)Y(t)d\Lambda(t)$, and $\int_0^{\tau}h_2(t)Y(t)dt$. Various Donsker results  in \citet[ch.19]{vdv1998} can be immediately applied here.  See the proofs of their Donsker properties in \citet[ch. 3.7]{kate2014}. Based on our assumptions as well as our definition of $H$  in  \eqref{function: h},  these classes are also all uniformly bounded. Applying the preservation theorem of Donsker \citep{van2000weak},  $\{\psi_{\theta, \Lambda,h}, (\theta, \Lambda) \in \mathbb{L}_0 \mbox{ and } h \in H \}$  is a Donsker class  with a finite envelop function. Since we have shown $\psi_{\theta, \Lambda, h}$ forms a Donsker class over the entire parameter space, 
condition \circled{2} on the GC property immediately holds by Slutsky's lemma. Now we verify condition \circled{4}. Since $P\psi_{\theta, \Lambda}$ is already a linear function of $\theta$ and $\Lambda$, it is straightforward to derive the Fr\'echet derivative $\dot{\Psi}_0$ of $P\psi_{\theta, \Lambda}$ at $(\theta_0, \Lambda_0)$ simply based on the definition of Fr\'echet derivative.  Written in a partition form, $\dot{\Psi}_0 = \left(\begin{array}{cc}
\dot{\Psi}_{11} &  \dot{\Psi}_{12}\\
 \dot{\Psi}_{21}&\dot{\Psi}_{22}
 \end{array} \right)$  
\begin{equation}\label{eqn:psi_dot}
\begin{aligned}
 \dot{\Psi}_{11}(\theta-\theta_0)h_1=& -h_1^TP\int_0^{\tau} Y(t)ZZ^T dt(\theta-\theta_0)\\
\dot{\Psi}_{12}(\Lambda-\Lambda_0)h_1=&  - h_1^TP \int_0^{\tau} Z Y(t)d(\Lambda- \Lambda_0)(t)\\
  \dot{\Psi}_{21}(\theta-\theta_0)h_2=& -P\int_0^{\tau} h_2(t)Y(t)Z^Tdt(\theta-\theta_0)\\
%lin \Theta  \mapsto l^{\infty}(H_2):    & (\theta-\theta_0)h_2 \mapsto
\dot{\Psi}_{22}(\Lambda-\Lambda_0)h_2= &-P\int_0^{\tau} h_2(t)Y(t)d(\Lambda-\Lambda_0)(t),
%lin \mathbb{A} \mapsto l^{\infty}(H_2)(\Lambda-\Lambda_0)h_2\mapsto
\end{aligned} 
\end{equation}

and $\dot{\Psi}_0$ is a map from $\mathbb{L}_0$ to $l^\infty(H)$. 
The proof for continuous invertibility of $\dot{\Psi}_0$ in a partitioned form can be  reduced to ascertaining the two operators $\dot{\Psi}_{22}$ and $\dot{V} = \dot{\Psi}_{11} - \dot{ \Psi}_{12}\dot{\Psi}_{22}^{-1}\dot{\Psi}_{21}$ have continuous inverses \citep[ch. 25]{vdv1998}. We are able to write both of their inverses  explicitly as below.
 Let ${\bf R}(\dot{ \Psi}_{22})$ denote the range of $\dot{ \Psi}_{22}$. Then 
\begin{equation}\label{psi22i}
(\dot{ \Psi}_{22}^{-1}\eta)h_2 = -\int_0^\tau \frac{h_2(t)}{PY(t)}d\eta(t).
\end{equation}
where $\eta \in {\bf R}(\dot{ \Psi}_{22})$. By assumption,  $PY(\tau)$ is bounded  below from 0. Thus,  $ \dot{ \Psi}_{22}^{-1} $ is continuous. Following the definition of $\dot{V}$, it is not difficult to derive that  $\dot{V}(\theta-\theta_0)h_1 = -h_1^TA(\theta-\theta_0)$ where $A$ is given in \eqref{eqn:A}. Since we assume $A$ is positive definite in A. \ref{assum: A-posdef}, its inverse $A^{-1}$ exists and so does $\dot{V}^{-1}$. $\dot{V}^{-1}$ as  a linear map between  Euclidean spaces is automatically continuous. Therefore, the inverse of $\dot{\Psi}_0$  exist and it is continuous. Finally, to verify condition \circled{5}, we can first show $\psi_{\theta,\Lambda,h} \rightarrow \psi_{\theta_0,\Lambda_0,h}$ pointwise and uniformly in $h$. Then, by dominated convergence theorem condition \circled{5} will hold. Interested readers are referred to the detailed proofs of these five conditions in \citet[ch. 3.7]{kate2014}.
%{\color{red}{I plan to write a separate article on additive hazards model with two-phase sampling alone with more details in theoretical steps and  simulation/validation results, so I keep this section brief. Is it OK to the journal? What do you think how detail I should give for the AH model example?}}

\subsection*{Step 4: Compute the Asymptotic Variances}
After conditions \circled{1} - \circled{5} are verified, applying corollaries \ref{thm: consist-srs}, \ref{thm: consist-2ph} and \ref{thm: consist-cali} we prove the consistency of $(\hat{\theta},\hat{\Lambda})$, $(\hat{\theta}^*,\hat{\Lambda}^*)$, and $(\hat{\theta}^{c}, \hat{\Lambda}^{c})$ immediately.  Applying theorems \ref{thm:gauss:rs},\ref{thm:gauss:2ph}, \ref{thm:gauss:cali}, we obtain the limiting  distributions of $\dot{\Psi}_0\sqrt {N}
\left(\begin{array}{c}
	\hat{\theta}-\theta_0 \\
	\hat{\Lambda}-\Lambda_0 \end{array}
\right)$, $\dot{\Psi}\sqrt {N}
\left(\begin{array}{c}
	\hat{\Lambda}^*-\Lambda_0\\ 
	\hat{\theta}^*-\theta_0 \end{array}
\right)$, and $\dot{\Psi}\sqrt {N}
\left(\begin{array}{c}
	\hat{\Lambda}^c-\Lambda_0\\ 
	\hat{\theta}^c-\theta_0 \end{array}
\right)$. 

Based on these joint limiting distributions,  the limiting distributions of  the estimators for the individual parameters of interest $\theta$, $\Lambda$, and $\Lambda(t|Z=z)$ are easy to obtain in a systematic way. The approach is as follows. We first find $h'= (h_1', h_2')$ such that $h' \in H$ and $\dot{\Psi}_0\sqrt{N}\left(\begin{array}{c}
\hat{\theta} - \theta_0 \\
\hat{\Lambda} - \Lambda_0 \end{array}
\right)h'$ is equal to the desired statistic $\sqrt{N}(\hat{\theta} - \theta_0)$,  $\sqrt{N}(\hat{\Lambda} - \theta_0)$, or $\sqrt{N}(\hat{\Lambda}(t|z) - \theta_0(t|z))$. Because the limiting distribution results in theorems \ref{thm:gauss:rs}, \ref{thm:gauss:2ph}, and \ref{thm:gauss:cali} hold for any $h \in H$, there is no need to prove asymptotic normality and derive the asymptotic variances for each individual estimator again. We only need to plug different $h'$s into the asymptotic variance formula in \eqref{var:srs} to retrieve the corresponding results.   Since $h'$ is determined solely by the parameter of interest, as data change to  two-phase VPS, we don't need to search for new $h'$ for two-phase VPS estimators or the calibrated version of them any more.  The only places that needs to be changed are the asymptotic theorems we invoke and the asymptotic variance formula we use,  i.e., \eqref{var:2ph} and \eqref{var:cali} instead.  The difference between deriving the model-based and robust variances lies only in the final step of the computation of the asymptotic variances.  The former one compute them under the assumption of \eqref{model:LY} and the latter not. They share the exact same procedure before this final step.

By this means, 9 statistics, with both model-based and robust variances,  from a combination of 3 parameters of interest and 3 data types can be derived very fast. Since this procedure is very standardized for all the statistics of our interest,  below we only show our development for $\hat{\theta}$, $\hat{\theta}^*$, and $\hat{\theta}^c$ as an illustration. We refer  the interested readers to \citet[ch. 3.10 \& ch. 4.3]{kate2014} for the development of $\hat{\Lambda}$, $\hat{\Lambda}^*$, $\hat{\Lambda}^c$, $\hat{\Lambda}(t|Z=z)$, $\hat{\Lambda}^*(t|Z=z)$, and $\hat{\Lambda}^c(t|Z=z)$. 

%The formula to calculate their asymptotic variances have also been given in equations \eqref{var:srs}, \eqref{var:2ph}, and \eqref{var:cali} respectively. 

%\subsection{limiting distribution of $\sqrt{N}(\hat{\theta}-\theta_0)$}\label{subsec:theta:srs}

\subsection*{Step A: find $h'$} 
In section \ref{sec:estimator},  we consider $\theta$ as an element in $l^{\infty}(H_1)$ identified by $\theta h_1=h_1^T \theta$. Hence, to find the limiting distribution of $\sqrt{N}(\hat{\theta}-\theta_0)$ is equivalent to finding
the asymptotic distribution of  $h_1^T\sqrt{N}(\hat{\theta}-\theta_0)$ for any $h_1 \in H_1$. Let
\begin{equation}\label{eqn:h:theta}
h' = (h_1',h_2') = \left(-A^{-1}h_1,  \frac{(A^{-1}h_1)^T PY(t)Z}{PY(t)} \right), t\in [0, \tau] 
\end{equation}
where $A$ is given in  \eqref{eqn:A}. Plugging $h'$ into the LHS of theorem \ref{thm:gauss:rs} yields
\begin{align*}\label{eqn:theta:left}
 &  \sqrt{N}\dot{\Psi}_{11}(\hat{\theta}-\theta_0)h_1'+\sqrt{N}\dot{\Psi}_{21}(\hat{\theta}-\theta_0)h_2'+ \sqrt{N}\dot{\Psi}_{12}(\hat{\Lambda}-\Lambda_0)h_1'+ \sqrt{N}\dot{\Psi}_{22} (\hat{\Lambda}-\Lambda_0)h_2'\nonumber\\
= & \sqrt{N}h_1^TA^{-1} \int_0^{\tau}\left[P\{Y(t)ZZ^T\}-\sqrt{N}\frac{P\{ Y(t)Z\}}{PY(t)}P\{Y(t)Z^T\}\right]dt(\hat{\theta}-\theta_0)+0\nonumber\\
=&\sqrt{N} h_1^TA^{-1}\int_0^{\tau}P\{Z-\frac{PZY(t)}{PY(t)}\}^{\otimes 2}Y(t) dt(\hat{\theta}-\theta_0)\nonumber\\
= & \sqrt{N}h_1^TA^{-1}A(\hat{\theta}-\theta_0)\nonumber\\
= & h_1^T\sqrt{N}(\hat{\theta}-\theta_0).
\end{align*}
$h'$ is motivated by aiming to cancel out $(\hat{\Lambda}- \Lambda_0)$ in $\dot{\Psi}_0\sqrt{N}\left(\begin{array}{c}
\hat{\theta} - \theta_0 \\
\hat{\Lambda} - \Lambda_0 \end{array}
\right)h'.$

\subsection*{Step B: show $h' \in H$}
 Because $h_1$ is  bounded and $A^{-1}$ exist, $\|-A^{-1}h_1\|_E$ is bounded and $h_1'$ still belong to $H_1$ in  \eqref{function: h}. %{\color{red}{ $<$- is this correct? what's the requirement for a matrix to be bounded?}}
 Next, we show $h_2' = \frac{(A^{-1}h_1)^T PY(t)Z}{PY(t)}  \in H_2$. By assumption $Z$ and $1/PY(\tau)$ are bounded, and $h_2'$ is uniformly bounded. Let $G$ denote the distribution of $Z$. Then for each coordinate of $Z_k, k = 1, 2, \dots, p$, 
\begin{align*}
PY(t)Z_k= & \int_{\mathbb{R}^p} E[1(T\geq t)|z]z_kdG(z)\\
      = & \int_{\mathbb{R}^p} S(t|z)z_kdG(z)\\
       = & \int_{\mathbb{R}^p, z_k>0} S(t|z)z_kdG(z) +\int_{\mathbb{R}^p, z_k <0} S(t|z)z_kdG(z).
\end{align*}
The first two term of the above display are monotonic in $t$, so $PY(t)Z_k$ as the summation of monotonic functions is of bounded variation.  $1/PY(t)$ is a bounded monotonic function of $t$, thus it is of $BV[0, \tau]$ as well. As the product of two bounded functions of bounded variation,  $\frac{PY(t)Z}{PY(t)}$ is of $BV[0, \tau]$. Hence $h'_2 \in H_2$. We have verified $h'\in H $.   

\subsection*{Step C: calculate the model-based and robust asymptotic variances}
Based on \eqref{eqn:normality_rs} and \eqref{var:srs},
$h'^T\sqrt{N}(\hat{\theta}-\theta_0)$ converges to a normal distribution with mean 0 and variance $P\psi^2_{\theta_0,\Lambda_0,h'}$. 
 Let $M(t)= N(t) - \Lambda_0(t) - Z^T\theta_0t$. Then 
 \[
 \psi_{\theta_0,\Lambda_0, h'}=  h_1'^T\psi_{1,\theta_0,\Lambda_0} + \psi_{2,\theta_0,\Lambda_0, h_2'} =  h_1'^T \int_{0}^{\tau}ZdM(t) +\int_{0}^{\tau}h_2'(t)dM(t)
 \]
and
\begin{align*}
 P\psi_{\theta_0,\Lambda_0,h'}^2
= & P\{h_1'^T \int_{0}^{\tau}ZdM(t) +\int_{0}^{\tau}h_2'(t)dM(t)\}^2\\
= & P\left\{-h_1^TA^{-1}\int_{0}^{\tau}ZdM(t) +\int_{0}^{\tau}\frac{h_1^TA^{-1}P Y(t)Z}{PY(t)}dM(t)\right\}^2 \\
=& h_1^TA^{-1}P\left[\int_{0}^{\tau}\left\{Z-\frac{PY(t)Z}{PY(t)}\right\}dM(t) \right]^{\otimes{2}}A^{-1}h_1.
\end{align*}
%On the RHS of \eqref{eqn:gauss:srs} we have 
%\begin{align}\label{eqn:theta:donsker}
% \mathbb{G}_N \left(\begin{array}{c}
%\psi_{1,\theta_0,\Lambda_0}\\
%\psi_{2,\theta_0,\Lambda_0}
%\end{array}\right)h'
% = &  \mathbb{G}_N (\psi_{1,\theta_0,\Lambda_0}h_1'+\psi_{2,\theta_0,\Lambda_0}h_2')\\
%  = & \mathbb{G}_N\left[-h_1^TA^{-1}\int_{0}^{\tau}ZdM(t) +\int_{0}^{\tau}\frac{h_1^TA^{-1}P Y(t)Z}{PY(t)}dM(t)\right]\nonumber\\
%  = & -\mathbb{G}_N h_1^TA^{-1}\int_{0}^{\tau}\left[Z-\frac{PY(t)Z}{PY(t)}\right] dM(t)
%\end{align}
%where 
%
%Thus 
%\begin{equation}\label{eqn:theta}
% \sqrt{N}(\hat{\theta}-\theta_0)h_1 =  -\mathbb{G}_N h_1^TA^{-1}\int_{0}^{\tau}\left[Z-\frac{PY(t)Z}{PY(t)}\right] dM(t).
%\end{equation}
Thus, 
$\sqrt{N}(\hat{\theta}-\theta_0)$ converges in distribution to a p-variate Gaussian distribution with mean zero and  variance  
\begin{equation*}\label{rvar:ah:theta}
 rVar_A\left\{ \sqrt{N}(\hat{\theta}-\theta_0)\right\} =  A^{-1}P\left[\int_{0}^{\tau}\left\{Z-\frac{PY(t)Z}{PY(t)}\right\}dM(t) \right]^{\otimes{2}}A^{-1}.
\end{equation*}
We denote this variance by $rVar$ to indicate this is the robust variance since its computation is carried out under the assumption of $P \in \mathcal{P}$ rather than  $P=P_{\theta_0,\Lambda_0}$  in \eqref{model:LY}. $M(t)$ is also not a martingale unless $P = P_{\theta_0,\Lambda_0}$.

When the AH model \eqref{model:LY} holds, $M(t)$ as the difference of a counting process $N(t)$ and an integrated intensity process $\int_0^{\cdot}Y(\cdot)d\Lambda(\cdot|Z=z)$ is a martingale. In consequence, $ \int_{0}^{\tau}\left\{Z-\frac{PY(t)Z}{PY(t)}\right\} Y(t)dM(t)$ is a  martingale integral with a predictable and locally bounded integrand. It follows from  the principle of martingale integral variance \citep[theorem 2.4.2]{FH} that 
\begin{align}\label{eqn:B}
  B \equiv & P\left[\int_{0}^{\tau}\left\{Z-\frac{PY(t)Z}{PY(t)}\right\} Y(t)dM(t) \right]^{\otimes{2}}\nonumber\\
=& Var\left[\int_{0}^{\tau}\left\{Z-\frac{PY(t)Z}{PY(t)}\right\}dM(t)\right] \nonumber\\
= & P\int_0^{\tau}\left\{Z-\frac{PY(t)Z}{PY(t)}\right\}^{\otimes 2}Y(t)\left\{d\Lambda_0(t)+Z^T\theta_0 dt\right\}
\end{align}
As a result, the model-based asymptotic variance
\begin{equation}\label{var:ah:theta}
  Var_A\left\{\sqrt{N}(\hat{\theta}-\theta_0)\right\}=  A^{-1}BA^{-1}.
\end{equation}
In conclusion,  when the AH model holds, $\sqrt{N}(\hat{\theta}-\theta_0)$  converges in distribution to a p-variate Gaussian distribution with mean zero and covariance matrix $A^{-1}BA^{-1}$. This result agrees with \citet{lin1994}. 

%\begin{align} 
 %Var_A\left\{ \sqrt{N}(\hat{\theta}-\theta_0)\right\} = & \nonumber A^{-1}P\left[\int_{0}^{\tau}\left[Z-\frac{P\{Y(t)Z\}}{PY(t)}\right]dM(t)\int_{0}^{\tau}\left[Z-\frac{P\{Y(t)Z\}}{PY(t)}\right]^TdM(t)\right]A^{-1} \nonumber\\
  %=& A^{-1} Var\left(\int_{0}^{\tau}\left[Z-\frac{P\{Y(t)Z\}}{PY(t)}\right]dM(t)\right)A^{-1}\nonumber\\
  %=& A^{-1}B A^{-1}.
%\end{align}

% $h_1^T\sqrt{N}(\hat{\theta}-\theta_0)$, we see that the technique is to find a pair of $(h_1',h_2') \in H$ such that 
%\[
%\left(\begin{array}{cc}
%\dot{ \Psi}_{11}& \dot{ \Psi}_{12}\\
%\dot{ \Psi}_{21}&\dot{\Psi}_{22} \end{array}
%\right) \sqrt {N}
%\left(\begin{array}{c}
%\hat{\theta}-\theta_0 \\
%\hat{\Lambda}-\Lambda_0 \end{array}
%\right)(h_1',h_2')= h_1^T\sqrt{N}(\hat{\theta}-\theta_0)
%\]
%After finding a suitable pair of $(h_1',h_2')$, we can apply results  \eqref{eqn:ah:gauss} and \eqref{var:ah} to estimate other statistics based on the joint distribution of $\sqrt{N}(\hat{\theta}-\theta_0, \hat{\theta}-\theta_0\}$. By adopting this approach,

%\subsection{limiting Distribution of $\sqrt{N}(\hat{\theta}^*-\theta_0)$ and $\sqrt{N}(\hat{\theta}^c-\theta_0)$}\label{subsec:theta:2ph}
Next, applying $h'$ in \eqref{eqn:h:theta} to theorem \ref{thm:gauss:2ph} and the result \eqref{var:2ph}, we obtain
$\sqrt{N}(\hat{\theta}^*-\theta_0)$  converges in distribution to a p-variate Gaussian distribution with mean zero and  variance  matrix
\begin{multline*}
 rVar_A\left\{ \sqrt{N}(\hat{\theta}^*-\theta_0)\right\} =  A^{-1}P\left\{\int_{0}^{\tau}\left[Z-\frac{PY(t)Z}{PY(t)}\right]dM(t) \right\}^{\otimes{2}}A^{-1}\\
 +Q\left[\frac{1-\pi_0(V)}{\pi_0(V)}A^{-1}\left[\int_{0}^{\tau}\left\{Z-\frac{PY(t)Z}{PY(t)}\right\}dM(t) \right]^{\otimes{2}}A^{-1}\right].
\end{multline*}
In view of \eqref{var:ah:theta},  when model \eqref{model:LY} holds the first term in the above display is replaced by the model-based variance for the RS estimator. Hence, the model-based asymptotic variance for $\hat{\theta}^*$ is
\begin{multline*}
  Var_A\left\{\sqrt{N}(\hat{\theta}^{*}-\theta_0)\right\}=  A^{-1}BA^{-1}\\
   + Q\left[\frac{1-\pi_0(V)}{\pi_0(V)}A^{-1}\left[\int_{0}^{\tau}\left\{Z-\frac{PY(t)Z}{PY(t)}\right\}dM(t) \right]^{\otimes{2}}A^{-1}\right].
\end{multline*}
Applying $h'$  in \eqref{eqn:h:theta} to  theorem \ref{thm:gauss:cali} and the the variance formula in \eqref{var:cali},  we obtain the robust asymptotic  variance  for the calibrated estimator $\hat{\theta}^c$ is
\begin{multline*}
rVar_A\left\{ \sqrt{N}(\hat{\theta}^c-\theta_0)\right\} = A^{-1}P\left\{\int_{0}^{\tau}\left[Z-\frac{PY(t)Z}{PY(t)}\right]dM(t) \right\}^{\otimes{2}}A^{-1}\\
+Q\left[\frac{1-\pi_0(V)}{\pi_0(V)}\left(A^{-1}\int_{0}^{\tau}\left\{Z-\frac{PY(t)Z}{PY(t)}\right\}dM\right.\right.\\ \left.\left.-Q\left[A^{-1}\int_{0}^{\tau}\left\{Z-\frac{PY(t)Z}{PY(t)}\right\}dM\tilde{V}^T\right](Q\tilde{V}\tilde{V}^T)^{-1}\tilde{V}\right)^{\otimes 2}\right]
\end{multline*}
and the model-based asymptotic variance   is
\begin{multline*}\label{var:theta:cali}
Var_A\left\{\sqrt{N}(\hat{\theta}^{c}-\theta_0)\right\}=  A^{-1}BA^{-1}\\   +Q\left[\frac{1-\pi_0(V)}{\pi_0(V)}\left(A^{-1}\int_{0}^{\tau}\left\{Z-\frac{PY(t)Z}{PY(t)}\right\}dM\right.\right.\\
\left.\left.-Q\left[A^{-1}\int_{0}^{\tau}\left\{Z-\frac{PY(t)Z}{PY(t)}\right\}dM\tilde{V}^T\right](Q\tilde{V}\tilde{V}^T)^{-1}\tilde{V}\right)^{\otimes 2}\right]. 
\end{multline*}
In summary, we systematically derived three estimators from three data types for $\theta$, each with two types of variances.  Using the same technique, one can derive estimators for $\Lambda_0$  with different data types and model assumptions by setting  \begin{equation*}\label{eqn:h:Lambda}
h' = (h_1',h_2') = \left(A^{-1}Dh_2,  -\frac{h_2(t)}{PY(t)}-\frac{(A^{-1}Dh_2)^T PY(t)Z}{PY(t)} \right),  t\in [0, \tau] 
\end{equation*} 
where  $D$ is defined by
\begin{equation*}\label{eqn:Dh2}
Dh_2 = \int_0^{\tau}h_2(t)\frac{PY(t)Z}{PY(t)}dt. 
\end{equation*}
 One can also derive the estimators for  $\Lambda(s|Z=z)$ with different data types and model assumptions   by setting \begin{equation*}\label{eqn:h:pred}
h'= (h_1',h_2') =  \left(A^{-1}\{D(s) - zs\},    -\frac{1(t \leq s)}{PY(t)}-\frac{\{D(s) - zs\}^T A^{-1}PY(t)Z}{PY(t)} \right)
\end{equation*} 
where $D(s) = Dh_2$ with $h_2 = 1(t \leq s)$.

\section{Discussion}

	The Z-estimation system was proposed to solve related estimation problems systematically. %for example, estimating the same parameter in a model with datasets from different sampling methods.  Although it is not necessary to use this system for all the estimation problems,  intermediate results from related problems are easily shared and reused if they are all developed in a single system. 
	In the following, we first summarize this system. 
	
\subsection{Summary}\label{sec:summary}
 The Z-estimation system divides an estimation procedure into four steps: 
	\begin{enumerate}
		\item [Step 1.] Define the parameter.  The true parameter $\alpha_0$ is defined as a functional that is usually motivated from a scientific question. This functional  is  in-explicitly  defined by a zero-valued equation $P\psi_{\alpha} = 0$ where $P$ is the true model. We allow this equation to be infinite-dimensional, so the parameter studied in this system can be infinite-dimensional too.  
		\item [Step 2.] Construct estimators. All the estimators for $\alpha_0$ are obtained by solving some EE, which are constructed in the form of  the average of  i.i.d. functions of $\psi_{\alpha}$ or a modification of it.   EE can be finite or infinite-dimensional depending on $\psi_{\alpha}$.
		\item [Step 3.]  Ascertain conditions \circled{1} - \circled{5} on $\psi_{\alpha}$. These conditions include \circled{1} $\alpha_0$ is a unique and well-separated solution to $P\psi_{\alpha} = 0$; \circled{2} $\psi_{\alpha}$ forms a GC class over the entire parameter space; \circled{3} $\psi_{\alpha}$ forms a Donsker class over the neighborhood of $\alpha_0$; \circled{4} The Fr\'echet derivative of $P\psi_{\alpha}$ and its continuous inverse exist; \circled{5} $\psi_{\alpha}$ converges to $\psi_{\alpha_0}$  in quadratic mean as $\alpha$ converges to $\alpha_0$. 
		\item [Step 4.] Calculate the asymptotic variances of these estimators. Select the appropriate consistency corollaries, asymptotic theorems, and asymptotic variances formula from section \ref{sec:tools} based on the data types and then assemble the ascertained five conditions to obtain the asymptotic results.  If $P = P_\alpha$,  model-based variances are computed by using this additional information,   otherwise robust variances are computed. 
\end{enumerate}
By defining parameters as a functional associated with the true model $P$, a proposed estimator is still interpretable when an assumed model $P_{\alpha}$ does not hold.  Such a treatment also helps unify the model-based and model-free estimation procedures. 
We allow the parameter to be Euclidean, non-Euclidean, or the combination of both. In consequence, this system can be  used not only for inference but also prediction if it is based on some semiparametric models as shown  in our  AH model example.

The three data types we studied for step 2 should be treated as a precursor for more interesting data types. For example, auxiliary information can come from a different data source rather than the phase I  of a two-phase study. %{\color{red}{Gary, are you aware other data problems I can mention here? I believe there is a lot of room for more research here. }} 
Consequently, the tool sets  in step 4 could be expanded when a new data type is encountered. The Z-estimation system was designed to ease this  expansion  as we exemplified in the development of our theoretical tools. Comparing our procedures in developing the asymptotic results for $\hat{\alpha} ^*$  to those for $\hat{\alpha}$, we adapted components \circled{2}\circled{3}\circled{5} to the new data types that incorporated auxiliary information while keeping other components unchanged. Integrating these small results together, the asymptotic results for $\hat{\alpha} ^*$ were obtained clearly and quickly.  We anticipate new assembly and theoretical tools could be built similarly for new data types.

\subsection{Limitation}

Methods of  weighted estimating equation and calibration are probably  the simplest  but may not be the most efficient estimation methods. This drawback and some methods for improvement  were discussed in \citet{breslow2007}.  In  certain scenarios, semiparmetric efficient estimators do exist. For example,  maximum likelihood or profile likelihood estimators are efficient and feasible if a semiparametric model can be partitioned into  parametric and non-parametric parts,  provided phase I data are discrete. \citet{lawless1999semiparametric, breslow2003,robins1995} also proposed optimal estimators for  the conditional  mean model when two-phase sampling is considered as a missing data problem with missing by design. We consider efficiency gain, however,  more from a perspective of using data intelligently during the data collection step  than data analytics techniques performed on a given dataset. For example, we improve efficiency  by selecting  the more  informative subjects  for expensive covariates measurements based on inexpensive information or collecting relevant auxiliary information. Our approach could be especially useful if technology makes the collection of inexpensive data  easy and fast, which is exactly what is happening now.   

This system also requires estimators to have $\sqrt{n}$ rate of convergence. When parameters are not estimable at $\sqrt{n}$-rate such as inference problems with interval  censoring data \citep{vdv1998} or mixture models \citep{bickel1993efficient},  this estimation system does not apply.	Alternative Z-theorems are available, which allow the convergence of nuisance parameters at a different rate with case-cohort data  \citep{nan2012}.
 %{\color{red}{are there new papers on this topic that I should be aware and cite here ? I see this is a major limitation of this approach. }}.

In this system, we only consider i.i.d. VPS to construct the phase II subsample in this system, so that our  observations are all independent. This restriction simplifies the asymptotic analysis, but the independent selection mechanism also causes variation in the size of the phase II subsamples within each stratum.  As a result, the variance of the estimator is slightly larger than the estimator that would be obtained if a fixed size were used to draw subsamples. However,  selecting the phase II subsamples with a fixed size generates correlation among them. Some theoretical tools were developed to handle this type of data \citep{breslow2007, saegusa2013}. Alternatively, the Z-estimation system is still applicable if we adapt this problem to a two-phase VPS problem with some calibration constraint. If we treat  the sample size for each phase II stratum as the auxiliary information, then the calibration constraint is that the  observed sampling fractions are equal to the ratios between the fixed phase II sample sizes and the phase I sample sizes for all the strata.  If we set the weight for each subject in EE \eqref{EE:calibrated_estimators} to be the inverse of the sampling fractions and use the indicator of a subject's stratum membership as the calibration variable, it can be shown that the resulting estimator and its variance derived from the Z-estimation system yield the same results as other methods. It is interesting to study whether a similar treatment could be considered for other types of correlated data. If so,  this system will be useful for broader data problems.
% since the constraints on sample size for each stratum was kept in the estimating procedure as the additional EE to be solved for calibrated estimators. 
%We sacrifice generality for simplicity by this restriction. 
\subsection{Systematic}
We take a ``decompose and assemble" approach  in three places within this system. The whole  estimation procedure was first divided into  four steps as summarized in section \ref{sec:summary}. Within step 3,  the asymptotic study of estimators is then further decomposed into five components. Among these five components, a similar ``decompose and assemble"   approach is also seen in verifying the GC and Donsker properties in  conditions \circled{2} and \circled{3} as illustrated in multiple proofs in sections \ref{sec:tools} and \ref{sec:example}. A studied function is first decomposed into several small basic functions. Then preservation theorems of GC and Donsker ensure that once these basic functions are shown to be GC and Donsker, these two properties are preserved for the function that integrates these basic functions by summation or multiplication under some minor conditions. %Several solutions to complex problems were taken into consideration while deciding the structure of this system. 

The four steps of our estimation procedure follow a natural thought process of solving a problem. On the other hand, they can also be viewed as almost  four independent parts that require different skills: step 1 on understanding and formulating a scientific problem, step 2 on data handling,  step 3 on asymptotic analysis, and step 4 on implementation.   In the past few decades, asymptotic analysis is often the main focus of an estimation problem  and a single well-trained statistician can usually conduct all the steps,  but today's problem requires large collaboration among professionals with diverse skills. We organize the estimation procedure to these four parts so that each part can be studied separately and thoroughly when they become complex. For example, in the past a  scientific problem  might had been already well formulated in previous literatures and the focus of the research was in methodology.  Now, however, digitalization has transformed each industry. New data are collected and expected to shed insights on a problem on which we never used data to study. A person who may not be familiar with asymptotic analysis techniques but has a profound domain knowledge can participate in step 1, suggesting  relevant and meaningful target parameters to estimate. Similarly, many data experts' strength lies in step 2 given today's big and complex data, and  software engineers' skills could be very useful  for creating a user-friendly, robust, and scalable API of estimators at step 4.   Therefore, we hope this system could be also used as a collaboration platform for domain experts, data experts, mathematicians, and software engineers. Their specialties will ensure a newly developed estimator is relevant, unbiased, reliable, and easy to be used.

%Before  data structure is often simple, data volume is not big, and  scientific problems may have been already well formulated in previous literature. Now, however, digitalization has transformed each industry. A large amount of data, with a complexed structure and from multiple sources are being collected, which requires data experts. To create a relevant new statistical metric or model to shed new insights on a scientific or business problem requires domain knowledge.  For a simple task, Step 4 may be carried out by a statistician in a R software, but if we would like a newly developed estimator to be widely used,  this step relies on the support of  software engineers to make the implementation robust and scalable. For example, how can we make a statistic calculated in real-time with real-time data that may be accessed by a large amount of people at one time? We purposely organize the estimation procedure to these four parts so each part has a room to be studied thoroughly. We hope this system could be also used as a collaboration platform from  domain experts, data experts, mathematician, and software engineers. Their specialties will ensure a new estimator is relevant to the problem, unbiased, reliable, and easy to be used.

Our estimating procedure can become a system because we chose to built it upon Z-estimation.
%We are able to have a systematic study of our proposed estimators  because  Z-estimation is by nature systematic.  Z-estimation started from  \citep{huber1967}. The old classic approach studied asymptotics of Z-estimators by the Taylor expansion on $\Psi_N(\alpha)$. Huber's Z-theorem instead separated the contribution to  $\Psi_N(\alpha)$ into  a deterministic part $\Psi(\alpha)$ and a stochastic remainder (see \citet{pollard1984} and \citet{pollard1985new}  for details). 	As a result, the conditions for asymptotics were decomposed into two components:   analytical conditions on $\Psi$ and stochastic conditions on the deviation of $\Psi$ from $\Psi_N$.  We recognized this "decompose and assemble" approach and  further decomposed the asymptotic analysis into conditions \circled{1} $\sim$ \circled{5}. 
The old classic   approach studied asymptotics of Z-estimators by Taylor expansion on $\Psi_N(\alpha)$. In contrast, Huber's 1967 landmark paper proposed a new approach for proof by separating the contribution to  $\Psi_N(\alpha)$ into  a deterministic part $\Psi(\alpha)$, which is  the expectation of  $\Psi_N(\alpha)$, and a stochastic remainder (see \citet[section VII]{pollard1984} and \citet{pollard1985new}  for details). 	As a result, the
conditions for asymptotics were decomposed into two components:   analytical conditions on $\Psi$ and stochastic conditions on the deviation of $\Psi$ from $\Psi_N$. %By using this systematic approach, Huber was able to relax conditions on  the second and higher order derivatives of likelihood function for ML estimators. 
%This new approach was not widely appreciated until  
Later Pollard [1984, 1985] made the connection of Huber's work with the  empirical process theory, which provides tools to check the  stochastic conditions. 
Today we call Huber (1967) and Pollard (1985)'s new way of proving the asymptotic normality of Z-estimators Z-estimation.  We recognized Z-estimation's  ``decompose and assemble" approach and  further decomposed the asymptotic analysis into smaller and smaller tasks.

 In this system, because we only consider independent observations and EE in $l^{\infty}$ space, the stochastic conditions for showing consistency and asymptotic normality are transformed to some GC and Donsker conditions, for which many existing results can be immediately used. The  GC and Donsker properties are themselves very convenient notations  for transferring knowledge.  Most EE are composed of some common basic functions. Once a basic function's GC and Donsker property has been proved somewhere,  the results can be directly used for new models. For example, in our AH model development, the basic functions of the EE include  $Z$, $\Delta$,  $T$, $Z^T\theta$, $\Lambda(T)$, $h_2(T)$, $\int_0^\tau h_2(t)Y(t)d\Lambda(t)$, and $\int_0^{\tau}h_2(t)Y(t)dt$. These are components commonly seen in the EE for many survival models. After we verified GC and Donsker properties of them  in our AH model development, researchers can immediately combine them with other functions  for new survival models development.  As the collection of the GC or Donsker functions become bigger, more and more theoretical work on checking stochastic conditions will be saved.There is a possibility that applied statisticians in the future may not need to verify GC or Donsker conditions directly themselves. By this means, the powerful and abstract modern empirical process theory results will become more  accessible and widely used for model development.

 Compared to  the asymptotic analysis with the Taylor expansion technique,  this system brings some transparency to the asymptotic analysis  by breaking it into small components. Small tasks in small components are easy to follow and check. This transparency, as a result, will improve the quality of methodology research. Small components also facilitate users to identify existing results, which might be reused and adapted for solving new problems. This system makes such adaptation easy as illustrated in our example when extending the RS estimators to the two-phase VPS estimators. %The new data types could be adapted,  the function in each dimension of the  EE could be more complex, and the  new constraints on estimation could be added %Finally, if we consider each GC and Donsker class of a basic function as a basic Lego game piece and the preservation theorems as the principles a user can rely on to bond these pieces,  the development of a new statistic could potentially become a Lego game. t
%The procedure of proving and deriving becomes a procedure of assembling. There are numerous way to assemble and interesting statistics might be built creatively by users themselves. 
The current information age demands creation of a large amount of statistics fast. This Z-estimation system, taking a systematic approach to build a series of related statistics on a single platform,  will expedite the development of statistics. As a system, it may reveal more  benefits as we use it.

		\section*{Appendix A}
		
		\begin{proof}[Proof of Lemma \ref{lemma:unique_solution}]
			
			We first show $0$ is the unique solution to the equation 
			$Q \psi^{c}_{2,\gamma}(V,R) 
			= \frac{R}{\pi_0(V)}exp(-\gamma^T \tilde{V}) \tilde{V} -\tilde{V} = 0$. By assumption A.\ref{assum: MAR}, we see $0$ is one solution to $Q\psi^{c}_{2,\gamma}(V,R)=0$.  Suppose there exists another solution $\gamma_1 \in \Gamma$ such that $\gamma_1\neq 0$ and $Q\psi^{c}_{2,\gamma_1}(V,R)=0$. 
			Then by the mean value theorem  in multiple variables
			\begin{align*}
			\psi^{c}_{2,\gamma_1}(V,R)-\psi^{c}_{2,0}(V,R)
			=&\frac{R}{\pi_0(V)}\exp(-\gamma_1^T \tilde{V}) \tilde{V} -\frac{R}{\pi_0(V)}exp(-\vec{0}^T \tilde{V}) \tilde{V}\\
			=&\frac{R}{\pi_0(V)}\exp(-\gamma_v^{*T}\tilde{V})\tilde{V}\tilde{V}^T(\gamma_1-\vec{0})
			\end{align*}
			where $\gamma_v^{*}$ is on the line segment between 0 and $\gamma_1$. We use notation $\gamma_v^{* }$ to suggest this vector depends on $V$. As a result, 
			\[
			Q\exp(-\gamma^{*T}_v\tilde{V})\tilde{V}\tilde{V}^T(\gamma_1-\vec{0})=
			Q\left[\psi^{c}_{2,\gamma_1}(V,R)-\psi^{c}_{2,0}(V,R)\right]=0
			.\]
			Given that $Q \tilde{V}\tilde{V}^T$ is positive definite, $ \gamma_1 ^TQ \tilde{V}\tilde{V}^T\gamma_1 >0$ unless 
			$\gamma_1=0$. Since $\Gamma$ is a compact convex subset of $\mathbb{R}^q$ with $0$ as an interior point, $\gamma_v^{*}$ belongs to $\Gamma$ and is bounded. By assumption A.\ref{assum: vtilde-bounded}, $\tilde{V}$ is also bounded. Thus $\exp(-\gamma^{*T
			}_v\tilde{V}) $ is always positive, bounded and bounded from 0. Therefore, $ \gamma_1 ^TQ \exp(-\gamma^{*T}_v\tilde{V})\tilde{V}\tilde{V}^T\gamma_1 >0$ unless 
			$\gamma_1=0$. By contradiction, we prove $\gamma=0$ is the unique solution to $Q\psi^{c}_{2,\gamma}(V,R)=0$.

		   Because $\alpha_0$ is the unique solution to the equation $\psi_{\alpha}(X) = 0$  and when $\gamma=0$ by \eqref{eqn:connection}
			\[
			Q\psi_{1,\alpha,\gamma}^{c}(R,V,X)=Q\frac{R}{\pi_0(V)}\psi_{\alpha}(X)=\psi_{\alpha}(X) = \Psi(\alpha),
			\]
		  $(\alpha_0,0)$ is the unique solution to the equations $Q\psi_{1,\alpha,\gamma}^{c}(R,V,X)=0$ and $Q\psi_{2, \gamma}^{c}(V,X) = 0$.
		\end{proof}
		
		\begin{proof}[Proof of Lemma \ref{fgamma}]
		
			Because $\Gamma$ is convex, for every $\gamma_1,\gamma_2 \in \Gamma$ we can find a $\gamma_v^* \in \Gamma $ such that
			\[
			|f_{\gamma_1}(R,X,V)-f_{\gamma_2}(R,X,V)| = |\frac{R}{\pi_0(V)}\exp (-\gamma_v^{*T}\tilde{V})\tilde{V}^T(\gamma_1-\gamma_2)|_\mathbb{E}.
			\]
			By Cauchy-Schwarz inequality,
			\[
			|\frac{R}{\pi_0(V)}\exp (-\gamma_v^{*T}\tilde{V})\tilde{V}^T(\gamma_1-\gamma_2)| \leq \| \frac{R}{\pi_0(V)}\exp (-\gamma_v^{*T}\tilde{V})\tilde{V}^T\|_\mathbb{E}\|\gamma_1-\gamma_2\|_\mathbb{E}.
			\]
			Since $\Gamma$ is compact, thus it is bounded. Under the assumption of A.\ref{assum: bounded_below} and A.\ref{assum: vtilde-bounded}, 
			both $1/\pi_0(V)$ and $\tilde{V}$ are also bounded. 
			Hence, there exists a positive constant $C$ such that 
			\[
			|f_{\gamma_1}(R,X,V)-f_{\gamma_2}(R,X,V)| \leq C\|\gamma_1-\gamma_2\|_\mathbb{E}.
			\]
			Furthermore, $\{ f_{\gamma}, \gamma \in \Gamma\}$ is uniformly bounded. As a result, $\{f_{\gamma},\gamma \in \Gamma \}$ is uniformly Lipschitz.
		\end{proof}
		
		\section*{Appendix B}	
		\begin{proof}[Verification of Fr\'echet derivative in Theorem \ref{thm:gauss:cali}]
			Let $\Psi_1^{c} (\alpha, \gamma) =  Q\psi^{c}_{1,\alpha,\gamma}(R,V,X)$ and $\Psi_2^{c} (\gamma) =  Q\psi^{c}_{2, \gamma}(R,V)$.
			We first verify $ \dot{\Psi}^{c}_{11}= \dot{\Psi}_0$ and $\dot{\Psi}^{c}_{21}= 0$.  By condition \circled{4}  $\dot{\Psi}_0$ as a continuous and linear map satisfies
			\[
			\|P\psi_{\alpha}(X)-P\psi_{\alpha_0}(X)- \dot{\Psi}_0(\alpha -\alpha_0)\|=o(\|\alpha-\alpha_0\|) \quad \mbox{as} \quad \|\alpha-\alpha_0\|\downarrow 0 .
			\]
			According to \eqref{eqn:connection}, $Q \psi^*_{\alpha}(R,V,X) = P \psi_{\alpha}(X)$. 
			Since  $\Psi_1^{c} (\alpha, \gamma_0) = Q\psi^*_{\alpha}(R,V,X)$,   $\Psi_1^{c} (\alpha, \gamma_0) =P \psi_{\alpha}(X)$.
			As a result, 
			\begin{align*}
			&\|\Psi_1^{c} (\alpha, \gamma)- \Psi_1^{c} (\alpha_0, \gamma)-\dot{\Psi}_0(\alpha-\alpha_0) \|\\
			= & \|\Psi_1^{c} (\alpha, \gamma)- \Psi_1^{c} (\alpha, \gamma_0)+\Psi_1^{c} (\alpha, \gamma_0) - \Psi_1^{c} (\alpha_0, \gamma_0)+ \Psi_1^{c} (\alpha_0, \gamma_0)-\Psi_1^{c} (\alpha_0, \gamma) -\dot{\Psi}_0(\alpha-\alpha_0) \| \\
			\leq &\|\Psi_1^{c} (\alpha, \gamma)- \Psi_1^{c} (\alpha, \gamma_0)\|+ \|\Psi_1^{c} (\alpha_0, \gamma_0)- \Psi_1^{c} (\alpha_0, \gamma)\|+\|\Psi_1^{c} (\alpha, \gamma_0)-\Psi_1^{c} (\alpha_0, \gamma_0)-\dot{\Psi}_0(\alpha-\alpha_0) \|\\
			= & 0+0+ \|P\psi_{\alpha}(X)-P\psi_{\alpha_0}(X)- \dot{\Psi}_0(\alpha -\alpha_0)\|\\
			=& o(\|\alpha-\alpha_0\|) \mbox{ at } (\alpha_0, \gamma_0) \mbox{ as } \|\alpha-\alpha_0\|\downarrow 0.
			\end{align*}
		 Therefore, $\dot{\Psi}_0$ is the Fr\'echet derivative of $\Psi_1^{c}$ with respect to $\alpha$ at $(\alpha_0,\gamma_0)$.
			
			 Because $\Psi_2^{c}(\gamma)$ does not involve $\alpha$, the Fr\'echet derivative of $\Psi_w^{c}$ with respect to $\alpha$ is $0$. Thus, $\dot{\Psi}^{c}_{21}= 0$.

			Next, we verify $ \dot{\Psi}^{c}_{12}=-Q\psi_{\alpha_0}\tilde{V} ^T$ and $ \dot{\Psi}^{c}_{22}=-Q\tilde{V}\tilde{V}^T$. For each $h \in \mathcal{H}$, we have 
			\begin{align*}
			\Psi_1^{c} (\alpha, \gamma)h & =Q\psi_{1, \alpha,\gamma,h}^{c}=
			Q\left[\frac{R}{\pi_0(V)}\exp(-\gamma^T\tilde{V})\psi_{\alpha,h }(X)\right] \\
			\Psi_2^{c} (\alpha, \gamma)& =Q\psi_{2, \gamma}^{c}=
			Q\left[\frac{R}{\pi_0(V)}\exp(-\gamma^T\tilde{V})\tilde{V}-\tilde{V}\right] \in \mathbb{R}^q.
			\end{align*}
			By assumption,  $\tilde{V}$ and $\gamma$ are bounded and under condition \circled{2}, $\{\psi_{\alpha,h},\alpha \in \mathbb{L}_0, h \in H\} $ has
			integrable envelope function. It follows from the dominated convergence theorem and 
			rules of differentiation, the gradients
			
			\begin{align*}
			\triangledown \dot{\Psi}_1^{c}(\alpha_0,\gamma_0)\cdot(\gamma-\gamma_0)h
			=& \triangledown Q\left[\frac{R}{\pi_0(V)}\exp(-\gamma^T\tilde{V})\psi_{\alpha,h}(X)\right]|_{\alpha_0, \gamma_0} \cdot (\gamma-\gamma_0)\\
			= & Q\left[ \triangledown \frac{R}{\pi_0(V)}\exp(-\gamma^T\tilde{V})\psi_{\alpha,h}(X)\right]|_{\alpha_0, \gamma_0} \cdot (\gamma-\gamma_0)\\
			= & -Q\left[ \psi_{\alpha_0,h}(X)\tilde{V}^T\right] )\cdot(\gamma-\gamma_0)
			\end{align*}
			
			and 
			\begin{align*}
			\triangledown \dot{\Psi}_2^{c}(\alpha_0,0)\cdot(\gamma-\gamma_0)
			=& \triangledown Q\left[\frac{R}{\pi_0(V)}\exp(-\gamma^T\tilde{V})\tilde{V}-\tilde{V}\right]|_{\alpha_0, \gamma_0} \cdot(\gamma-\gamma_0)\\
			= & Q\triangledown \left[\frac{R}{\pi_0(V)}\exp(-\gamma^T\tilde{V})\tilde{V}-\tilde{V}\right]|_{\alpha_0, \gamma_0} \cdot (\gamma-\gamma_0)\\
			= & -Q[ \tilde{V}\tilde{V}^T] \cdot(\gamma-\gamma_0).
			\end{align*}
		%	{\color{red}{Does the above notation make sense? If I remember it correctly, this step is based on Norm's college note!}}
			Because on a Euclidean space, Fr\'echet derivative agrees with differential, we have 
			\[
			\dot{\Psi}_{22}^{c}(\gamma -\gamma_0)= -Q[ \tilde{V}\tilde{V}^T] (\gamma-\gamma_0).
			\]
			For each $h \in H$, by Taylor expansion and dominated convergence theorem
			\begin{align*}
			\Psi_1^{c} (\alpha, \gamma)h= & \Psi_1^{c} (\alpha, \gamma_0)h- Q\left[ \psi_{\alpha,h}\tilde{V}^T\right] (\gamma-\gamma_0)\\
			& \qquad +1/2 (\gamma-\gamma_0)^TQ[\tilde{V}\tilde{V}^T \psi_{\alpha,h}\exp(-\gamma_v^T\tilde{V})](\gamma-\gamma_0)
			\end{align*}
			where $\gamma_v$ is some line segment joining 0 and $\gamma$.
			%Given each $\gamma \in \Gamma$, We denote the map $h \mapsto Q\left[ \psi_{\alpha_0,h}(X)\tilde{V}^T\right] (\gamma-\gamma_0)$ by 
			%$Q\left[ \psi_{\alpha_0}\tilde{V}^T\right] (\gamma-\gamma_0)$.
			Then
			\begin{align*}
			&\|\Psi_1^{c} (\alpha, \gamma)- \Psi_1^{c} (\alpha, \gamma_0)+Q\left[ \psi_{\alpha_0}\tilde{V}^T\right] (\gamma-\gamma_0) \|\\
			=&\sup_{h \in H}|\Psi_1^{c} (\alpha, \gamma)h- \Psi_1^{c} (\alpha, \gamma_0)h +Q\left[ \psi_{\alpha_0,h}(X)\tilde{V}^T\right] (\gamma-\gamma_0) |\\
			= & (\gamma-\gamma_0)^T\sup_{h \in H}\left|1/2 Q[\tilde{V}\tilde{V}^T \psi_{\alpha_0,h}\exp(-\gamma_v^T\tilde{V})]\right|(\gamma-\gamma_0)\\
			=& o(\| \gamma-\gamma_0\|) \mbox{ as } \| \gamma-\gamma_0\| \downarrow 0.
			\end{align*}
			
			Therefore, $-Q\psi_{\alpha_0}\tilde{V} ^T$ is the Fr\'echet derivative of $\Psi_1^{c}$ with respect to $\gamma$ at $(\alpha_0,\gamma_0)$ and $-Q\tilde{V}\tilde{V}^T$ is the Fr\'echet derivative of $\Psi_2^{c}$ with respect to $\gamma$
			at $(\alpha_0,\gamma_0)$. 
			
			In summary, the Fr\'echet derivative of the map $(\alpha,\gamma) \mapsto \Psi^{c} (\alpha,\gamma)$ at $(\alpha_0,\gamma_0)$ is
			\[\dot{\Psi}^{c}_0 =\left(\begin{array}{cc}
			\dot{\Psi}^{c}_{11} & \dot{\Psi}^{c}_{12}\\
			\dot{\Psi}^{c}_{21} & \dot{\Psi}^{c}_{22}
			\end{array}\right)
			= \left(\begin{array}{cc}
			\dot{\Psi}_0 & 
			-Q\psi_{\alpha_0}(X)\tilde{V}^T\\
			0 &
			-Q\tilde{V}\tilde{V}^T
			\end{array}\right).
			\]
		\end{proof}

 \bibliographystyle{imsart-nameyear}
 \bibliography{draft}

\end{document}